\documentclass{article}
\usepackage{amsfonts}
\usepackage{amsmath}

\def\limfunc#1{\mathop{\rm #1}}

\begin{document}

\title{Accompanying document to\\
\textquotedblleft Point Estimation with
Exponentially Tilted Empirical Likelihood\textquotedblright }
\author{Susanne M. Schennach}
\date{First draft: December 12, 2002\\
This version: December 7, 2005}
\maketitle

This manuscript provide the calculational details involved in the proof of
the equivalence between the $O\left( n^{-2}\right) $ variance of the ETEL
and EL estimator that are omitted from the appendix of \textquotedblleft
Point Estimation with Exponentially Tilted Empirical
Likelihood\textquotedblright\ by Susanne M.\ Schennach.

Whenever possible, the following calculations were verified using the
symbolic calculations\ capabilities of Maple.

\section{Definitions}

Let 
\begin{eqnarray*}
\beta &=&\left( \tau ,\kappa ^{\prime },\lambda ^{\prime },\theta ^{\prime
}\right) ^{\prime } \\
\beta ^{\ast } &=&\left( 1,\mathbf{0}^{\prime },\mathbf{0}^{\prime },\theta
^{\ast \prime }\right) ^{\prime }
\end{eqnarray*}

\begin{eqnarray*}
\dot{g} &=&g\left( x_{i},\theta \right) \\
\dot{G} &=&\frac{\partial g\left( x_{i},\theta \right) }{\partial \theta
^{\prime }}
\end{eqnarray*}%
To simplify the notation, the dependence of the above quantities on $i$ and $%
\theta $ is implicit.%
\begin{eqnarray*}
G &=&E\left[ G\left( x_{i},\theta ^{\ast }\right) \right] \\
\Omega &=&E\left[ g\left( x_{i},\theta ^{\ast }\right) g\left( x_{i},\theta
^{\ast }\right) ^{\prime }\right]
\end{eqnarray*}%
\begin{eqnarray*}
\bar{g} &=&\left[ n^{-1/2}\sum_{i}\dot{g}\right] _{\beta =\beta ^{\ast }} \\
\bar{G} &=&\left[ n^{-1/2}\sum_{i}\left( \dot{G}-G\right) \right] _{\beta
=\beta ^{\ast }} \\
\bar{\Omega} &=&\overline{gg^{\prime }}=\left[ n^{-1/2}\sum_{i}\left( \dot{g}%
\dot{g}^{\prime }-\Omega \right) \right] _{\beta =\beta ^{\ast }}
\end{eqnarray*}%
Convention for transpose: $G_{\cdot h}^{\prime }\equiv \left( G_{\cdot
h}\right) ^{\prime }$

Expectations evaluated at $\beta =\beta ^{\ast }$ will be denoted by $%
E^{\ast }\left[ \ldots \right] $. For instance, $E^{\ast }\left[ \dot{G}%
\right] =G$ and $E^{\ast }\left[ \dot{g}\dot{g}^{\prime }\right] =\Omega $.

\subsection{Moment conditions for ETEL}

In ETEL, $\hat{\beta}^{ETEL}$ solves $n^{-1}\sum_{i}\dot{\phi}^{ETEL}=0$,
where

\[
\dot{\phi}^{ETEL}=\left[ 
\begin{array}{c}
\dot{\tau}-\tau \\ 
\dot{\tau}\dot{g} \\ 
\left( \tau -\dot{\tau}\right) \dot{g}+\dot{\tau}\dot{g}\dot{g}^{\prime
}\kappa \\ 
\dot{\tau}\dot{G}^{\prime }\kappa +\dot{\tau}\dot{G}^{\prime }\lambda \dot{g}%
^{\prime }\kappa -\dot{\tau}\dot{G}^{\prime }\lambda +\tau \dot{G}^{\prime
}\lambda%
\end{array}%
\right]
\]%
and where%
\[
\dot{\tau}=\exp \left( \lambda ^{\prime }\dot{g}\right) .
\]%
Define%
\begin{eqnarray*}
\Phi _{l,j}^{ETEL} &=&E^{\ast }\left[ \frac{\partial \dot{\phi}_{l}^{ETEL}}{%
\partial \beta _{j}}\right] \\
\Phi _{l,jk}^{ETEL} &=&E^{\ast }\left[ \frac{\partial ^{2}\dot{\phi}%
_{l}^{ETEL}}{\partial \beta _{j}\partial \beta _{k}}\right] \\
\Phi _{l,jkh}^{ETEL} &=&E^{\ast }\left[ \frac{\partial ^{3}\dot{\phi}%
_{l}^{ETEL}}{\partial \beta _{j}\partial \beta _{k}\partial \beta _{h}}%
\right]
\end{eqnarray*}%
\begin{eqnarray*}
\bar{\Phi}_{l}^{ETEL} &=&\left[ n^{-1/2}\sum_{i}\dot{\phi}_{l}^{ETEL}\right]
_{\beta =\beta ^{\ast }} \\
\bar{\Phi}_{l,j}^{ETEL} &=&\left[ n^{-1/2}\sum_{i}\left( \frac{\partial \dot{%
\phi}_{l}^{ETEL}}{\partial \beta _{j}}-\Phi _{l,j}^{ETEL}\right) \right]
_{\beta =\beta ^{\ast }} \\
\bar{\Phi}_{l,jk}^{ETEL} &=&\left[ n^{-1/2}\sum_{i}\left( \frac{\partial ^{2}%
\dot{\phi}_{l}^{ETEL}}{\partial \beta _{j}\partial \beta _{k}}-\Phi
_{l,jk}^{ETEL}\right) \right] _{\beta =\beta ^{\ast }} \\
\bar{\Phi}_{l,jkh}^{ETEL} &=&\left[ n^{-1/2}\sum_{i}\left( \frac{\partial
^{3}\dot{\phi}_{l}^{ETEL}}{\partial \beta _{j}\partial \beta _{k}\partial
\beta _{h}}-\Phi _{l,jkh}^{ETEL}\right) \right] _{\beta =\beta ^{\ast }}
\end{eqnarray*}%
\begin{eqnarray*}
\bar{\Phi}_{\cdot ,\cdot }^{ETEL} &=&\left[ 
\begin{array}{cccc}
-1 & 0 & \bar{g}^{\prime } & 0 \\ 
0 & 0 & \bar{\Omega} & \bar{G} \\ 
\bar{g} & \bar{\Omega} & -\bar{\Omega} & 0 \\ 
0 & \bar{G}^{\prime } & 0 & 0%
\end{array}%
\right] _{\beta =\beta ^{\ast }} \\
\Phi _{\cdot ,\cdot }^{ETEL} &=&\left[ 
\begin{array}{cccc}
-1 & 0 & 0 & 0 \\ 
0 & 0 & \Omega & G \\ 
0 & \Omega & -\Omega & 0 \\ 
0 & G^{\prime } & 0 & 0%
\end{array}%
\right]
\end{eqnarray*}%
where the dot replacing the subscripts denote all the elements of the matrix.

\subsection{Moment conditions for EL}

The moment conditions for the EL parameters $\left( \hat{\kappa}^{\prime },%
\hat{\theta}^{\prime }\right) ^{\prime }$ are 
\[
n^{-1}\sum_{i}\left[ 
\begin{array}{c}
\dot{\varepsilon}\dot{g} \\ 
\dot{\varepsilon}\dot{G}^{\prime }\kappa%
\end{array}%
\right] =0
\]%
where $\dot{\varepsilon}=\left( 1-\kappa ^{\prime }\dot{g}\right) ^{-1}$ and 
$\kappa $ is the Lagrange multiplier of the moment contraints, which has
been relabelled $\kappa $ to simplify the comparison with ETEL. Furthermore,
to again simplify the comparison with ETEL, we augment this vector by $%
1+\dim \kappa $ additional moment conditions and introduce the same number
of additional parameters $\left( \tau ,\lambda \right) $ where $\tau \in 
\mathbb{R}$ and $\lambda \in \mathbb{R}^{\dim \kappa }$: 
\[
n^{-1}\sum_{i}\left[ 
\begin{array}{c}
\dot{\tau}-\tau \\ 
\dot{\tau}\dot{g} \\ 
\dot{\varepsilon}\dot{g} \\ 
\dot{\varepsilon}\dot{G}^{\prime }\kappa%
\end{array}%
\right] =0
\]%
where $\dot{\tau}=\exp \left( \lambda ^{\prime }\dot{g}\right) $. The
additional moment conditions merely define the values of the new parameters $%
\left( \hat{\tau},\hat{\lambda}\right) $ and do not change the values of $%
\left( \hat{\kappa}^{\prime },\hat{\theta}^{\prime }\right) $. Indeed,
whenever $\left( \hat{\kappa}^{\prime },\hat{\theta}^{\prime }\right)
^{\prime }$ are such that the bottom two subvectors are zero, one can always
find a value of $\left( \hat{\tau},\hat{\lambda}\right) $ that will make the
top two subvectors vanish as well. Since the origin is in the convex hull of 
$\left\{ g\left( x_{i},\hat{\theta}\right) \right\} _{i=1}^{n}$ w.p.a. 1,
there exists $\hat{\lambda}$ such that $n^{-1}\sum_{i}\dot{\tau}\dot{g}=0$
w.p.a. 1. Then, we can just set $\hat{\tau}=n^{-1}\sum_{i}\dot{\tau}$.

Finally, since just-identified GMM is invariant under linear transformations
of the vector of moment conditions, the moment conditions for EL can
equivalently be written as $n^{-1}\sum_{i}\dot{\phi}^{EL}=0$, where%
\[
\dot{\phi}^{EL}=\left[ 
\begin{array}{c}
\dot{\tau}-\tau \\ 
\dot{\tau}\dot{g} \\ 
\dot{\varepsilon}\dot{g}-\dot{\tau}\dot{g} \\ 
\dot{\varepsilon}\dot{G}^{\prime }\kappa%
\end{array}%
\right]
\]%
This particular version of the EL moment conditions will drastically
simplify our calculations, due to the fact that the matrices of first
derivatives for ETEL and EL become nearly identical:%
\begin{eqnarray*}
\bar{\Phi}_{\cdot ,\cdot }^{EL} &=&\left[ 
\begin{array}{cccc}
-1 & 0 & \bar{g}^{\prime } & 0 \\ 
0 & 0 & \bar{\Omega} & \bar{G} \\ 
0 & \bar{\Omega} & -\bar{\Omega} & 0 \\ 
0 & \bar{G}^{\prime } & 0 & 0%
\end{array}%
\right] \\
\Phi _{\cdot ,\cdot }^{EL} &=&\left[ 
\begin{array}{cccc}
-1 & 0 & 0 & 0 \\ 
0 & 0 & \Omega & G \\ 
0 & \Omega & -\Omega & 0 \\ 
0 & G^{\prime } & 0 & 0%
\end{array}%
\right]
\end{eqnarray*}%
Note that 
\[
\Phi _{\cdot ,\cdot }^{ETEL}=\Phi _{\cdot ,\cdot }^{EL}\equiv \Phi _{\cdot
,\cdot }
\]%
while $\bar{\Phi}_{\cdot ,\cdot }^{ETEL}$ and $\bar{\Phi}_{\cdot ,\cdot
}^{EL}$ differ by a single element.

\subsection{Other definitions and conventions}

\begin{eqnarray*}
\Psi _{\cdot ,j}^{ETEL} &=&-\left( \Phi _{\cdot ,\cdot }\right) ^{-1}\Phi
_{\cdot ,j}^{ETEL} \\
\Psi _{\cdot ,jk}^{ETEL} &=&-\left( \Phi _{\cdot ,\cdot }\right) ^{-1}\Phi
_{\cdot ,jk}^{ETEL} \\
\Psi _{\cdot ,jkh}^{ETEL} &=&-\left( \Phi _{\cdot ,\cdot }\right) ^{-1}\Phi
_{\cdot ,jkh}^{ETEL}
\end{eqnarray*}%
\begin{eqnarray*}
\bar{\Psi}_{\cdot }^{ETEL} &=&-\left( \Phi _{\cdot ,\cdot }\right) ^{-1}\bar{%
\Phi}_{\cdot }^{ETEL} \\
\bar{\Psi}_{\cdot ,j}^{ETEL} &=&-\left( \Phi _{\cdot ,\cdot }\right) ^{-1}%
\bar{\Phi}_{\cdot ,j}^{ETEL} \\
\bar{\Psi}_{\cdot ,jk}^{ETEL} &=&-\left( \Phi _{\cdot ,\cdot }\right) ^{-1}%
\bar{\Phi}_{\cdot ,jk}^{ETEL}
\end{eqnarray*}%
and similarly for EL.

Again, a dot replacing a given subscript denotes a vector of all the values
taken for all values of that subscript.

e.g. $\Psi _{\cdot ,j}^{ETEL}=\left[ 
\begin{array}{c}
\Psi _{1,j}^{ETEL} \\ 
\vdots \\ 
\Psi _{1+2\dim \lambda +\dim \theta ,j}^{ETEL}%
\end{array}%
\right] $, $\Psi _{l,\cdot }^{ETEL}=\left[ 
\begin{array}{ccc}
\Psi _{l,1}^{ETEL} & \cdots & \Psi _{l,1+2\dim \lambda +\dim \theta }^{ETEL}%
\end{array}%
\right] $.

Following Newey and Smith (2001), define 
\begin{eqnarray*}
P &=&\Omega ^{-1}-\Omega ^{-1}G\left( G^{\prime }\Omega ^{-1}G\right)
^{-1}G^{\prime }\Omega ^{-1} \\
H &=&\left( G^{\prime }\Omega ^{-1}G\right) ^{-1}G^{\prime }\Omega ^{-1} \\
\Sigma &=&\left( G^{\prime }\Omega ^{-1}G\right) ^{-1}
\end{eqnarray*}

Identities:%
\begin{eqnarray*}
PG &=&0 \\
P^{\prime } &=&P \\
P\Omega P &=&P \\
P\Omega H^{\prime } &=&0 \\
H\Omega H^{\prime } &=&\Sigma
\end{eqnarray*}%
Define%
\begin{eqnarray*}
l_{\tau } &=&0 \\
l_{\kappa } &=&1 \\
l_{\lambda } &=&1+\dim \kappa \\
l_{\theta } &=&1+2\dim \kappa
\end{eqnarray*}%
These symbols will be used to isolate subvectors and submatrices. For
instance, 
\[
\Phi _{h,l_{\kappa }+j}^{EL}=E\left[ \frac{\partial \dot{\phi}_{h}^{EL}}{%
\partial \kappa _{j}}\right]
\]%
with $j$ implicitly ranging from $1$ to $\dim \kappa $.

\section{Stochastic expansion}

The following merely rewrites the conclusion of Lemma A4 in Newey and Smith
(2001).using our notation.

\[
\left( \hat{\beta}-\beta ^{\ast }\right) =n^{-1/2}\bar{\Psi}_{\cdot }+n^{-1}%
\bar{Q}_{\cdot }+n^{-3/2}\bar{R}_{\cdot }+O_{p}\left( n^{-2}\right)
\]

where%
\begin{eqnarray*}
\bar{Q}_{l} &=&\bar{\Psi}_{l,j}\bar{\Psi}_{j}+\frac{1}{2}\Psi _{l,jk}\bar{%
\Psi}_{j}\bar{\Psi}_{k} \\
\bar{R}_{l} &=&\bar{\Psi}_{l,j}\bar{Q}_{j}+\Psi _{l,jk}\bar{Q}_{j}\bar{\Psi}%
_{k}+\frac{1}{2}\bar{\Psi}_{l,jk}\bar{\Psi}_{j}\bar{\Psi}_{k}+\frac{1}{6}%
\Psi _{l,jkh}\bar{\Psi}_{j}\bar{\Psi}_{k}\bar{\Psi}_{h}
\end{eqnarray*}%
using the convention that repeated indices are summed over, e.g. $\bar{\Psi}%
_{l,j}\bar{\Psi}_{j}\equiv \sum_{j}\bar{\Psi}_{l,j}\bar{\Psi}_{j}$.

(We also simplified Newey and Smith's result using the fact that $\Psi
_{l,jk}=\Psi _{l,kj}$).

The quantities associated with each estimator will distinguished by an
\textquotedblleft ETEL\textquotedblright\ or \textquotedblleft
EL\textquotedblright\ superscript.

The variance of $\hat{\beta}$ is:%
\[
\limfunc{Var}\left[ \left( \hat{\beta}-\beta ^{\ast }\right) \right] =n^{-1}%
\limfunc{Var}\left[ \bar{\Psi}_{\cdot }\right] +n^{-2}\limfunc{Var}\left[ 
\bar{Q}_{\cdot }\right] +n^{-2}\text{Covar}\left[ \bar{R}_{\cdot },\bar{\Psi}%
_{\cdot }^{\prime }\right] +n^{-2}\text{Covar}\left[ \bar{\Psi}_{\cdot },%
\bar{R}_{\cdot }^{\prime }\right] +o\left( n^{-2}\right)
\]

\section{Partitioned Inverse of $\Phi _{\cdot ,\cdot }$}

$\left[ 
\begin{array}{cccc}
-1 & 0 & 0 & 0 \\ 
0 & 0 & \Omega & G \\ 
0 & \Omega & -\Omega & 0 \\ 
0 & G^{\prime } & 0 & 0%
\end{array}%
\right] ^{-1}=\left[ 
\begin{array}{cc}
-1 & \left[ 
\begin{array}{ccc}
0 & 0 & 0%
\end{array}%
\right] \\ 
\left[ 
\begin{array}{c}
0 \\ 
0 \\ 
0%
\end{array}%
\right] & \left[ 
\begin{array}{ccc}
0 & \Omega & G \\ 
\Omega & -\Omega & 0 \\ 
G^{\prime } & 0 & 0%
\end{array}%
\right] ^{-1}%
\end{array}%
\right] $

where $\left[ 
\begin{array}{ccc}
0 & \Omega & G \\ 
\Omega & -\Omega & 0 \\ 
G^{\prime } & 0 & 0%
\end{array}%
\right] ^{-1}$ can be found using partitioned inverse formula:

$\left[ 
\begin{array}{cc}
A_{11} & A_{12} \\ 
A_{21} & A_{22}%
\end{array}%
\right] ^{-1}=\left[ 
\begin{array}{cc}
B_{11} & B_{12} \\ 
B_{21} & B_{22}%
\end{array}%
\right] $

$A_{11}^{-1}=\left[ 
\begin{array}{cc}
0 & \Omega \\ 
\Omega & -\Omega%
\end{array}%
\right] ^{-1}=\left[ 
\begin{array}{cc}
\Omega ^{-1} & \Omega ^{-1} \\ 
\Omega ^{-1} & 0%
\end{array}%
\right] $

$B_{22}=\left( A_{22}-A_{21}A_{11}^{-1}A_{12}\right) ^{-1}=\left( 0-\left[ 
\begin{array}{cc}
G^{\prime } & 0%
\end{array}%
\right] \left[ 
\begin{array}{cc}
\Omega ^{-1} & \Omega ^{-1} \\ 
\Omega ^{-1} & 0%
\end{array}%
\right] \left[ 
\begin{array}{c}
G \\ 
0%
\end{array}%
\right] \right) ^{-1}=-\left( G^{\prime }\Omega ^{-1}G\right) ^{-1}\equiv
-\Sigma $

$B_{11}=A_{11}^{-1}\left( I+A_{12}B_{22}A_{21}A_{11}^{-1}\right) =\left[ 
\begin{array}{cc}
\Omega ^{-1} & \Omega ^{-1} \\ 
\Omega ^{-1} & 0%
\end{array}%
\right] \left( I-\left[ 
\begin{array}{c}
G \\ 
0%
\end{array}%
\right] \Sigma \left[ 
\begin{array}{cc}
G^{\prime } & 0%
\end{array}%
\right] \left[ 
\begin{array}{cc}
\Omega ^{-1} & \Omega ^{-1} \\ 
\Omega ^{-1} & 0%
\end{array}%
\right] \right) =$

$=\left[ 
\begin{array}{cc}
\Omega ^{-1} & \Omega ^{-1} \\ 
\Omega ^{-1} & 0%
\end{array}%
\right] \left[ 
\begin{array}{cc}
I-G\Sigma G^{\prime }\Omega ^{-1} & -G\Sigma G^{\prime }\Omega ^{-1} \\ 
0 & I%
\end{array}%
\right] \allowbreak =\left[ 
\begin{array}{cc}
\Omega ^{-1}-\Omega ^{-1}G\Sigma G^{\prime }\Omega ^{-1} & \Omega
^{-1}-\Omega ^{-1}G\Sigma G^{\prime }\Omega ^{-1} \\ 
\Omega ^{-1}-\Omega ^{-1}G\Sigma G^{\prime }\Omega ^{-1} & -\Omega
^{-1}G\Sigma G^{\prime }\Omega ^{-1}%
\end{array}%
\right] $

$=\left[ 
\begin{array}{cc}
P & P \\ 
P & P-\Omega ^{-1}%
\end{array}%
\right] $

$B_{12}=-A_{11}^{-1}A_{12}B_{22}=-\left[ 
\begin{array}{cc}
\Omega ^{-1} & \Omega ^{-1} \\ 
\Omega ^{-1} & 0%
\end{array}%
\right] \left[ 
\begin{array}{c}
G \\ 
0%
\end{array}%
\right] \Sigma =-\left[ 
\begin{array}{c}
\Omega ^{-1}G\Sigma \\ 
\Omega ^{-1}G\Sigma%
\end{array}%
\right] =\left[ 
\begin{array}{c}
H^{\prime } \\ 
H^{\prime }%
\end{array}%
\right] $

$B_{21}=B_{12}^{\prime }$

$\left[ 
\begin{array}{cccc}
-1 & 0 & 0 & 0 \\ 
0 & 0 & \Omega & G \\ 
0 & \Omega & -\Omega & 0 \\ 
0 & G^{\prime } & 0 & 0%
\end{array}%
\right] ^{-1}=\left[ 
\begin{array}{cccc}
-1 & 0 & 0 & 0 \\ 
0 & P & P & H^{\prime } \\ 
0 & P & P-\Omega ^{-1} & H^{\prime } \\ 
0 & H & H & -\Sigma%
\end{array}%
\right] $

For notational convenience, we define%
\[
\Phi _{jk}^{-1}\equiv \left( \Phi _{\cdot ,\cdot }^{-1}\right) _{jk}
\]

\section{$\bar{\Psi}_{\cdot }$Term}

$\bar{\Psi}_{\cdot }^{ETEL}=\bar{\Psi}_{\cdot }^{EL}\equiv \bar{\Psi}_{\cdot
}=-\Phi _{\cdot ,\cdot }^{-1}\bar{\Phi}_{\cdot }=-\left[ 
\begin{array}{cccc}
-1 & 0 & 0 & 0 \\ 
0 & P & P & H^{\prime } \\ 
0 & P & P-\Omega ^{-1} & H^{\prime } \\ 
0 & H & H & -\Sigma%
\end{array}%
\right] \left[ 
\begin{array}{c}
0 \\ 
\bar{g} \\ 
0 \\ 
0%
\end{array}%
\right] =\left[ 
\begin{array}{c}
0 \\ 
-P\bar{g} \\ 
-P\bar{g} \\ 
-H\bar{g}%
\end{array}%
\right] $

$\limfunc{Var}\left( \bar{\Psi}_{\cdot }\right) =E\left[ \left[ 
\begin{array}{c}
0 \\ 
-P\bar{g} \\ 
-P\bar{g} \\ 
-H\bar{g}%
\end{array}%
\right] \left[ 
\begin{array}{cccc}
0 & -\bar{g}^{\prime }P & -\bar{g}^{\prime }P & -\bar{g}^{\prime }H^{\prime }%
\end{array}%
\right] \right] =\left[ 
\begin{array}{cccc}
0 & 0 & 0 & 0 \\ 
0 & P & P & 0 \\ 
0 & P & P & 0 \\ 
0 & 0 & 0 & \Sigma%
\end{array}%
\right] $

\section{$\bar{Q}_{\cdot }$ TERM}

\[
\bar{Q}_{l}=\bar{\Psi}_{l,j}\bar{\Psi}_{j}+\bar{\Psi}_{j}\Psi _{l,jk}\bar{%
\Psi}_{k}/2
\]

\subsection{Proof that $\bar{Q}_{l}^{ETEL}-\bar{Q}_{l}^{EL}=0$}

This result will imply that all the parameters $\left( \hat{\tau},\hat{\kappa%
}^{\prime },\hat{\lambda}^{\prime },\hat{\theta}\right) $ for ETEL and EL
differ by less than $O_{p}\left( n^{-1}\right) $. This may come as a
surprise since the Lagrange multipliers of EL and ET are known to differ by $%
O_{p}\left( n^{-1}\right) $. However, we have carefully defined the
parameter vector of each estimator so that the Lagrange multiplier of EL ($%
\hat{\kappa}$) corresponds to the auxiliary parameter $\hat{\kappa}$ of ETEL
instead of the Lagrange multiplier $\hat{\lambda}$ of ETEL.

\[
\left[ 
\begin{array}{r}
\hat{\tau}^{ETEL} \\ 
\text{Auxiliary parameter of ETEL }\rightarrow \hat{\kappa}^{ETEL} \\ 
\text{Lagrange multiplier of ETEL }\rightarrow \hat{\lambda}^{ETEL} \\ 
\hat{\theta}^{ETEL}%
\end{array}%
\right] -\left[ 
\begin{array}{l}
\hat{\tau}^{EL} \\ 
\hat{\kappa}^{EL}\text{ }\leftarrow \text{ Lagrange multiplier of EL} \\ 
\hat{\lambda}^{EL}\text{ }\leftarrow \text{ parameter added for convenience}
\\ 
\hat{\theta}^{EL}%
\end{array}%
\right]
\]%
Hence, our results do not imply that the Lagrange multiplier EL and ETEL are
equivalent up to $O_{p}\left( n^{-1}\right) $ but rather that the auxiliary
parameter $\hat{\kappa}^{ETEL}$ in fact plays the role of EL's Lagrange
multiplier in the ETEL estimator.

Since the influence function of EL and ETEL are the same, we have that

\[
\bar{Q}_{l}^{ETEL}-\bar{Q}_{l}^{EL}=\left( \bar{\Psi}_{l,j}^{ETEL}-\bar{\Psi}%
_{l,j}^{EL}\right) \bar{\Psi}_{j}+\bar{\Psi}_{j}\left( \Psi
_{l,jk}^{ETEL}-\Psi _{l,jk}^{EL}\right) \bar{\Psi}_{k}/2
\]

\subsubsection{Calculation of $\left( \bar{\Psi}_{l,j}^{ETEL}-\bar{\Psi}%
_{l,j}^{EL}\right) \bar{\Psi}_{j}$}

$\left( \bar{\Psi}_{l,j}^{ETEL}-\bar{\Psi}_{l,j}^{EL}\right) \bar{\Psi}%
_{j}=\Phi _{lh}^{-1}\left( \bar{\Phi}_{h,j}^{ETEL}-\bar{\Phi}%
_{h,j}^{EL}\right) \bar{\Psi}_{j}$, where

$\left( \bar{\Phi}_{\cdot ,\cdot }^{ETEL}-\bar{\Phi}_{\cdot ,\cdot
}^{EL}\right) \bar{\Psi}_{\cdot }=\left( \left[ 
\begin{array}{cccc}
-1 & 0 & \bar{g}^{\prime } & 0 \\ 
0 & 0 & \bar{g}\bar{g}^{\prime } & \bar{G} \\ 
\bar{g} & \bar{g}\bar{g}^{\prime } & -\bar{g}\bar{g}^{\prime } & 0 \\ 
0 & \bar{G}^{\prime } & 0 & 0%
\end{array}%
\right] -\left[ 
\begin{array}{cccc}
-1 & 0 & \bar{g}^{\prime } & 0 \\ 
0 & 0 & \bar{g}\bar{g}^{\prime } & \bar{G} \\ 
0 & \bar{g}\bar{g}^{\prime } & -\bar{g}\bar{g}^{\prime } & 0 \\ 
0 & \bar{G}^{\prime } & 0 & 0%
\end{array}%
\right] \right) \left[ 
\begin{array}{c}
0 \\ 
-P\bar{g} \\ 
-P\bar{g} \\ 
-H\bar{g}%
\end{array}%
\right] =$

$\qquad =\left[ 
\begin{array}{cccc}
0 & 0 & 0 & 0 \\ 
0 & 0 & 0 & 0 \\ 
\bar{g} & 0 & 0 & 0 \\ 
0 & 0 & 0 & 0%
\end{array}%
\right] \left[ 
\begin{array}{c}
0 \\ 
-P\bar{g} \\ 
-P\bar{g} \\ 
-H\bar{g}%
\end{array}%
\right] =\left[ 
\begin{array}{c}
0 \\ 
0 \\ 
0 \\ 
0%
\end{array}%
\right] $

Thus, $\left( \bar{\Psi}_{\cdot ,\cdot }^{ETEL}-\bar{\Psi}_{\cdot ,\cdot
}^{EL}\right) \bar{\Psi}_{\cdot }=\mathbf{0}.$

\subsubsection{Calculation of $\bar{\Psi}_{j}\left( \Psi _{l,jk}^{ETEL}-\Psi
_{l,jk}^{EL}\right) \bar{\Psi}_{k}/2$}

\label{secpsiquad}$\bar{\Psi}_{j}\left( \Psi _{l,jk}^{ETEL}-\Psi
_{l,jk}^{EL}\right) \bar{\Psi}_{k}=\Phi _{lh}^{-1}\bar{\Psi}_{j}\left( \Phi
_{h,jk}^{ETEL}-\Phi _{h,jk}^{EL}\right) \bar{\Psi}_{k}=\Phi _{lh}^{-1}\bar{%
\Psi}_{\cdot }^{\prime }E^{\ast }\left[ \frac{\partial ^{2}}{\partial \beta
\partial \beta ^{\prime }}\left( \phi _{h}^{ETEL}-\phi _{h}^{EL}\right) %
\right] \bar{\Psi}_{\cdot }$

where

$\left( \dot{\phi}^{ETEL}-\dot{\phi}^{EL}\right) =\left( \left[ 
\begin{array}{c}
\dot{\tau}-\tau \\ 
\dot{\tau}\dot{g} \\ 
\left( \tau -\dot{\tau}\right) \dot{g}+\dot{\tau}\dot{g}\dot{g}^{\prime
}\kappa \\ 
\dot{\tau}\dot{G}^{\prime }\kappa +\dot{\tau}\dot{G}^{\prime }\lambda \dot{g}%
^{\prime }\kappa -\dot{\tau}\dot{G}^{\prime }\lambda +\tau \dot{G}^{\prime
}\lambda%
\end{array}%
\right] -\left[ 
\begin{array}{c}
\dot{\tau}-\tau \\ 
\dot{\tau}\dot{g} \\ 
\dot{\varepsilon}\dot{g}-\dot{\tau}\dot{g} \\ 
\dot{\varepsilon}\dot{G}^{\prime }\kappa%
\end{array}%
\right] \right) $

$\qquad =\left[ 
\begin{array}{c}
0 \\ 
0 \\ 
\left( \tau -1\right) \dot{g}+\left( \dot{\tau}-\dot{\varepsilon}\right) 
\dot{g}\dot{g}^{\prime }\kappa \\ 
\left( \dot{\tau}-\dot{\varepsilon}\right) \dot{G}^{\prime }\kappa +\dot{\tau%
}\dot{G}^{\prime }\lambda \dot{g}^{\prime }\kappa -\left( \dot{\tau}-\tau
\right) \dot{G}^{\prime }\lambda%
\end{array}%
\right] $

where we used the identity $\dot{\varepsilon}\dot{g}=\frac{\dot{g}}{1-\dot{g}%
^{\prime }\kappa }=\dot{g}\left( 1+\frac{\dot{g}^{\prime }\kappa }{1-\dot{g}%
^{\prime }\kappa }\right) =\dot{g}\left( 1+\dot{\varepsilon}\dot{g}^{\prime
}\kappa \right) =\dot{g}+\dot{\varepsilon}\dot{g}\dot{g}^{\prime }\kappa $.

To calculate $\left( \Phi _{h,jk}^{ETEL}-\Phi _{h,jk}^{EL}\right) $:

$\left[ \frac{\partial ^{2}}{\partial \beta \partial \beta ^{\prime }}\left( 
\dot{\phi}_{l_{\lambda }+h}^{ETEL}-\dot{\phi}_{l_{\lambda }+h}^{EL}\right) %
\right] _{\beta =\beta ^{\ast }}=\left[ \frac{\partial ^{2}}{\partial \beta
\partial \beta ^{\prime }}\left( \left( \tau -1\right) \dot{g}_{h}+\left( 
\dot{\tau}-\dot{\varepsilon}\right) \dot{g}_{h}\dot{g}^{\prime }\kappa
\right) \right] _{\beta =\beta ^{\ast }}=$

$=\left[ 
\begin{array}{cccc}
0 & 0 & 0 & \dot{G}_{h\cdot } \\ 
0 & -2\dot{g}\dot{g}^{\prime }\dot{g}_{h} & \dot{g}_{h}\dot{g}\dot{g}%
^{\prime } & 0 \\ 
0 & \dot{g}_{h}\dot{g}\dot{g}^{\prime } & 0 & 0 \\ 
\dot{G}_{h\cdot }^{\prime } & 0 & 0 & 0%
\end{array}%
\right] _{\beta =\beta ^{\ast }}$

and

$\left[ \frac{\partial ^{2}}{\partial \beta \partial \beta ^{\prime }}\left( 
\dot{\phi}_{l_{\theta }+h}^{ETEL}-\dot{\phi}_{l_{\theta }+h}^{EL}\right) %
\right] _{\beta =\beta ^{\ast }}=\left[ \frac{\partial ^{2}}{\partial \beta
\partial \beta ^{\prime }}\left( \left( \dot{\tau}-\dot{\varepsilon}\right) 
\dot{G}_{\cdot h}^{\prime }\kappa +\dot{\tau}\dot{G}_{\cdot h}^{\prime
}\lambda \dot{g}^{\prime }\kappa +\left( \tau -\dot{\tau}\right) \dot{G}%
_{\cdot h}^{\prime }\lambda \right) \right] _{\beta =\beta ^{\ast }}$

$\qquad =\left[ 
\begin{array}{cccc}
0 & 0 & \dot{G}_{\cdot h}^{\prime } & 0 \\ 
0 & -\frac{\partial \left( \dot{g}\dot{g}^{\prime }\right) }{\partial \theta
_{h}} & \frac{\partial \left( \dot{g}\dot{g}^{\prime }\right) }{\partial
\theta _{h}} & 0 \\ 
\dot{G}_{\cdot h} & \frac{\partial \left( \dot{g}\dot{g}^{\prime }\right) }{%
\partial \theta _{h}} & -\frac{\partial \left( \dot{g}\dot{g}^{\prime
}\right) }{\partial \theta _{h}} & 0 \\ 
0 & 0 & 0 & 0%
\end{array}%
\right] _{\beta =\beta ^{\ast }}$

Then

$\bar{\Psi}_{j}\left( \Phi _{l_{\tau }+h,jk}^{ETEL}-\Phi _{l_{\tau
}+h,jk}^{EL}\right) \bar{\Psi}_{k}=0$

$\bar{\Psi}_{j}\left( \Phi _{l_{\kappa }+h,jk}^{ETEL}-\Phi _{l_{\kappa
}+h,jk}^{EL}\right) \bar{\Psi}_{k}=0$

$\bar{\Psi}_{j}\left( \Phi _{l_{\lambda }+h,jk}^{ETEL}-\Phi _{l_{\lambda
}+h,jk}^{EL}\right) \bar{\Psi}_{k}=$

$\qquad =\left[ 
\begin{array}{cccc}
0 & -\bar{g}^{\prime }P & -\bar{g}^{\prime }P & -\bar{g}^{\prime }H^{\prime }%
\end{array}%
\right] \left[ 
\begin{array}{cccc}
0 & 0 & 0 & G_{h\cdot } \\ 
0 & -2E^{\ast }\left[ \dot{g}\dot{g}^{\prime }\dot{g}_{h}\right] & E^{\ast } 
\left[ \dot{g}_{h}\dot{g}\dot{g}^{\prime }\right] & 0 \\ 
0 & E^{\ast }\left[ \dot{g}_{h}\dot{g}\dot{g}^{\prime }\right] & 0 & 0 \\ 
G_{h\cdot }^{\prime } & 0 & 0 & 0%
\end{array}%
\right] \left[ 
\begin{array}{c}
0 \\ 
-P\bar{g} \\ 
-P\bar{g} \\ 
-H\bar{g}%
\end{array}%
\right] =0$

\qquad where $G_{\cdot h}^{\prime }\equiv \left( G_{\cdot h}\right) ^{\prime
}$.

$\bar{\Psi}_{j}\left( \Phi _{l_{\theta }+h,jk}^{ETEL}-\Phi _{l_{\theta
}+h,jk}^{EL}\right) \bar{\Psi}_{k}=$

$\qquad =\left[ 
\begin{array}{cccc}
0 & -\bar{g}^{\prime }P & -\bar{g}^{\prime }P & -\bar{g}^{\prime }H^{\prime }%
\end{array}%
\right] 
\begin{array}{cccc}
0 & 0 & G_{\cdot h}^{\prime } & 0 \\ 
0 & -E^{\ast }\left[ \frac{\partial \left( \dot{g}\dot{g}^{\prime }\right) }{%
\partial \theta _{h}}\right] & E^{\ast }\left[ \frac{\partial \left( \dot{g}%
\dot{g}^{\prime }\right) }{\partial \theta _{h}}\right] & 0 \\ 
G_{\cdot h} & E^{\ast }\left[ \frac{\partial \left( \dot{g}\dot{g}^{\prime
}\right) }{\partial \theta _{h}}\right] & -E^{\ast }\left[ \frac{\partial
\left( \dot{g}\dot{g}^{\prime }\right) }{\partial \theta _{h}}\right] & 0 \\ 
0 & 0 & 0 & 0%
\end{array}%
\left[ 
\begin{array}{c}
0 \\ 
-P\bar{g} \\ 
-P\bar{g} \\ 
-H\bar{g}%
\end{array}%
\right] $

$\qquad =0$

Since $\bar{\Psi}_{j}\left( \Phi _{h,jk}^{ETEL}-\Phi _{h,jk}^{EL}\right) 
\bar{\Psi}_{k}=0$ for all $h$, it follows that $\bar{\Psi}_{j}\left( \Psi
_{l,jk}^{ETEL}-\Psi _{l,jk}^{EL}\right) \bar{\Psi}_{k}=0$ for all $l$.

Conclusion: $\bar{Q}_{l}^{ETEL}-\bar{Q}_{l}^{EL}=0$.

\subsection{Calculation of $\bar{Q}_{l}^{ETEL}=\bar{Q}_{l}^{EL}\equiv \bar{Q}%
_{l}$}

\[
\bar{Q}_{l}^{EL}=\bar{\Psi}_{l,j}^{EL}\bar{\Psi}_{j}+\bar{\Psi}_{j}\Psi
_{l,jk}^{EL}\bar{\Psi}_{k}/2
\]

\subsubsection{Calculation of $\bar{\Psi}_{l,j}^{EL}\bar{\Psi}_{j}$}

$\bar{\Psi}_{\cdot ,\cdot }^{EL}=-\left[ 
\begin{array}{cccc}
-1 & 0 & 0 & 0 \\ 
0 & P & P & H^{\prime } \\ 
0 & P & \left( P-\Omega ^{-1}\right) & H^{\prime } \\ 
0 & H & H & -\Sigma%
\end{array}%
\right] \left[ 
\begin{array}{cccc}
0 & 0 & \bar{g}^{\prime } & 0 \\ 
0 & 0 & \bar{\Omega} & \bar{G} \\ 
0 & \bar{\Omega} & -\bar{\Omega} & 0 \\ 
0 & \bar{G}^{\prime } & 0 & 0%
\end{array}%
\right] =-\left[ 
\begin{array}{cccc}
0 & 0 & -\bar{g}^{\prime } & 0 \\ 
0 & P\bar{\Omega}+H^{\prime }\bar{G}^{\prime } & 0 & P\bar{G} \\ 
0 & \left( P-\Omega ^{-1}\right) \bar{\Omega}+H^{\prime }\bar{G}^{\prime } & 
\Omega ^{-1}\bar{\Omega} & P\bar{G} \\ 
0 & H\bar{\Omega}-\Sigma \bar{G}^{\prime } & 0 & H\bar{G}%
\end{array}%
\right] \allowbreak $

$\bar{\Psi}_{\cdot ,\cdot }^{EL}\bar{\Psi}_{\cdot }=-\left[ 
\begin{array}{cccc}
0 & 0 & -\bar{g}^{\prime } & 0 \\ 
0 & P\bar{\Omega}+H^{\prime }\bar{G}^{\prime } & 0 & P\bar{G} \\ 
0 & P\bar{\Omega}-\Omega ^{-1}\bar{\Omega}+H^{\prime }\bar{G}^{\prime } & 
\Omega ^{-1}\bar{\Omega} & P\bar{G} \\ 
0 & H\bar{\Omega}-\Sigma \bar{G}^{\prime } & 0 & H\bar{G}%
\end{array}%
\right] \left[ 
\begin{array}{c}
0 \\ 
-P\bar{g} \\ 
-P\bar{g} \\ 
-H\bar{g}%
\end{array}%
\right] =\left[ 
\begin{array}{c}
-\bar{g}^{\prime }P\bar{g} \\ 
P\bar{\Omega}P\bar{g}+H^{\prime }\bar{G}^{\prime }P\bar{g}+P\bar{G}H\bar{g}
\\ 
P\bar{\Omega}P\bar{g}+H^{\prime }\bar{G}^{\prime }P\bar{g}+P\bar{G}H\bar{g}
\\ 
H\bar{\Omega}P\bar{g}-\Sigma \bar{G}^{\prime }P\bar{g}+H\bar{G}H\bar{g}%
\end{array}%
\right] $

\subsubsection{Calculation of $\bar{\Psi}_{j}\Psi _{l,jk}^{EL}\bar{\Psi}%
_{k}/2$}

$\bar{\Psi}_{j}\Psi _{l,jk}^{EL}\bar{\Psi}_{k}/2=-\Phi _{lh}^{-1}\bar{\Psi}%
_{j}\Phi _{h,jk}^{EL}\bar{\Psi}_{k}/2$ where

$\bar{\Psi}_{j}\Phi _{l_{\tau }+h,jk}^{EL}\bar{\Psi}_{k}=\bar{\Psi}%
_{j}E^{\ast }\left[ \frac{\partial ^{2}}{\partial \beta _{j}\partial \beta
_{k}}\dot{\phi}_{l_{\tau }+h}^{EL}\right] \bar{\Psi}_{k}=\bar{\Psi}_{\cdot
}^{\prime }E^{\ast }\left[ \frac{\partial ^{2}\left( \dot{\tau}-\tau \right) 
}{\partial \beta \partial \beta ^{\prime }}\right] \bar{\Psi}_{\cdot }=$

$\qquad =\left[ 
\begin{array}{cccc}
0 & -\bar{g}^{\prime }P & -\bar{g}^{\prime }P & -\bar{g}^{\prime }H^{\prime }%
\end{array}%
\right] \left[ 
\begin{array}{cccc}
0 & 0 & 0 & 0 \\ 
0 & 0 & 0 & 0 \\ 
0 & 0 & \Omega & G \\ 
0 & 0 & G^{\prime } & 0%
\end{array}%
\right] \left[ 
\begin{array}{c}
0 \\ 
-P\bar{g} \\ 
-P\bar{g} \\ 
-H\bar{g}%
\end{array}%
\right] =g^{\prime }\left( P\Omega P\right) g+g^{\prime }H^{\prime }\left(
G^{\prime }P\right) g+g^{\prime }\left( PG\right) Hg$

$\qquad =\bar{g}^{\prime }P\bar{g}+0+0$

$\bar{\Psi}_{j}\Phi _{l_{\kappa }+h,jk}^{EL}\bar{\Psi}_{k}=\bar{\Psi}%
_{j}E^{\ast }\left[ \frac{\partial ^{2}}{\partial \beta _{j}\partial \beta
_{k}}\dot{\phi}_{l_{\kappa }+h}^{EL}\right] \bar{\Psi}_{k}=\bar{\Psi}_{\cdot
}^{\prime }E^{\ast }\left[ \frac{\partial ^{2}\left( \dot{\tau}\dot{g}%
_{h}\right) }{\partial \beta \partial \beta ^{\prime }}\right] \bar{\Psi}%
_{\cdot }=$

$\qquad =\left[ 
\begin{array}{cccc}
0 & -\bar{g}^{\prime }P & -\bar{g}^{\prime }P & -\bar{g}^{\prime }H^{\prime }%
\end{array}%
\right] \left[ 
\begin{array}{cccc}
0 & 0 & 0 & 0 \\ 
0 & 0 & 0 & 0 \\ 
0 & 0 & E^{\ast }\left[ \dot{g}_{h}\dot{g}\dot{g}^{\prime }\right] & E^{\ast
}\left[ \frac{\partial }{\partial \theta ^{\prime }}\left( \dot{g}_{h}\dot{g}%
\right) \right] \\ 
0 & 0 & E^{\ast }\left[ \frac{\partial }{\partial \theta }\left( \dot{g}_{h}%
\dot{g}^{\prime }\right) \right] & E^{\ast }\left[ \frac{\partial ^{2}\dot{g}%
_{h}}{\partial \theta \partial \theta ^{\prime }}\right]%
\end{array}%
\right] \left[ 
\begin{array}{c}
0 \\ 
-P\bar{g} \\ 
-P\bar{g} \\ 
-H\bar{g}%
\end{array}%
\right] =$

$\qquad =\bar{g}^{\prime }PE^{\ast }\left[ \dot{g}_{h}\dot{g}\dot{g}^{\prime
}\right] P\bar{g}+\bar{g}^{\prime }PE^{\ast }\left[ \frac{\partial }{%
\partial \theta ^{\prime }}\left( \dot{g}_{h}\dot{g}\right) \right] H\bar{g}+%
\bar{g}^{\prime }H^{\prime }E^{\ast }\left[ \frac{\partial }{\partial \theta 
}\left( \dot{g}_{h}\dot{g}^{\prime }\right) \right] P\bar{g}+\bar{g}^{\prime
}H^{\prime }E^{\ast }\left[ \frac{\partial ^{2}\dot{g}_{h}}{\partial \theta
\partial \theta ^{\prime }}\right] H\bar{g}$

$\qquad =\left( \left( \bar{g}^{\prime }P\right) _{j}E^{\ast }\left[ \dot{g}%
\dot{g}_{j}\dot{g}^{\prime }\right] P\bar{g}\right) _{h}+\left( \left( \bar{g%
}^{\prime }P\right) _{j}E^{\ast }\left[ \frac{\partial }{\partial \theta
^{\prime }}\left( \dot{g}\dot{g}_{j}\right) \right] H\bar{g}\right)
_{h}+\left( \left( \bar{g}^{\prime }H^{\prime }\right) _{j}E^{\ast }\left[ 
\frac{\partial }{\partial \theta _{j}}\left( \dot{g}\dot{g}^{\prime }\right) %
\right] P\bar{g}\right) _{h}+$

$\qquad \qquad +\left( \left( \bar{g}^{\prime }H^{\prime }\right)
_{j}E^{\ast }\left[ \frac{\partial ^{2}\dot{g}}{\partial \theta _{j}\partial
\theta ^{\prime }}\right] H\bar{g}\right) _{h}$

$\bar{\Psi}_{j}\Phi _{l_{\lambda }+h,jk}^{EL}\bar{\Psi}_{k}=\bar{\Psi}%
_{j}E^{\ast }\left[ \frac{\partial ^{2}}{\partial \beta _{j}\partial \beta
_{k}}\dot{\phi}_{l_{\lambda }+h}^{EL}\right] _{\beta =\beta ^{\ast }}\bar{%
\Psi}_{k}=\bar{\Psi}_{\cdot }^{\prime }E^{\ast }\left[ \frac{\partial
^{2}\left( \dot{\varepsilon}\dot{g}_{h}-\dot{\tau}\dot{g}_{h}\right) }{%
\partial \beta \partial \beta ^{\prime }}\right] _{\beta =\beta ^{\ast }}%
\bar{\Psi}_{\cdot }=$

$\qquad =\left[ 
\begin{array}{cccc}
0 & -\bar{g}^{\prime }P & -\bar{g}^{\prime }P & -\bar{g}^{\prime }H^{\prime }%
\end{array}%
\right] \left[ 
\begin{array}{cccc}
0 & 0 & 0 & 0 \\ 
0 & 2E^{\ast }\left[ \dot{g}\dot{g}^{\prime }\dot{g}_{h}\right] & 0 & 
E^{\ast }\left[ \frac{\partial \left( \dot{g}\dot{g}_{h}\right) }{\partial
\theta ^{\prime }}\right] \\ 
0 & 0 & -E^{\ast }\left[ \dot{g}\dot{g}^{\prime }\dot{g}_{h}\right] & 
-E^{\ast }\left[ \frac{\partial \left( \dot{g}\dot{g}_{h}\right) }{\partial
\theta ^{\prime }}\right] \\ 
0 & E^{\ast }\left[ \frac{\partial \left( \dot{g}_{h}\dot{g}^{\prime
}\right) }{\partial \theta }\right] & -E^{\ast }\left[ \frac{\partial \left( 
\dot{g}_{h}\dot{g}^{\prime }\right) }{\partial \theta }\right] & 0%
\end{array}%
\right] \left[ 
\begin{array}{c}
0 \\ 
-P\bar{g} \\ 
-P\bar{g} \\ 
-H\bar{g}%
\end{array}%
\right] =$

$\qquad =\bar{g}^{\prime }PE^{\ast }\left[ \dot{g}\dot{g}^{\prime }\dot{g}%
_{h}\right] P\bar{g}=\left( \left( \bar{g}^{\prime }P\right) _{j}E^{\ast }%
\left[ \dot{g}_{j}\dot{g}\dot{g}^{\prime }\right] P\bar{g}\right) _{h}$

$\bar{\Psi}_{j}\Phi _{l_{\theta }+h,jk}^{EL}\bar{\Psi}_{k}=\bar{\Psi}%
_{j}E^{\ast }\left[ \frac{\partial ^{2}}{\partial \beta _{j}\partial \beta
_{k}}\dot{\phi}_{l_{\theta }+h}^{EL}\right] _{\beta =\beta ^{\ast }}\bar{\Psi%
}_{k}=\bar{\Psi}_{\cdot }^{\prime }E^{\ast }\left[ \frac{\partial ^{2}\left( 
\dot{\varepsilon}\dot{G}_{\cdot h}^{\prime }\kappa \right) }{\partial \beta
\partial \beta ^{\prime }}\right] _{\beta =\beta ^{\ast }}\bar{\Psi}_{\cdot
}=$

$\qquad =\left[ 
\begin{array}{cccc}
0 & -\bar{g}^{\prime }P & -\bar{g}^{\prime }P & -\bar{g}^{\prime }H^{\prime }%
\end{array}%
\right] \left[ 
\begin{array}{cccc}
0 & 0 & 0 & 0 \\ 
0 & E^{\ast }\left[ \frac{\partial \left( \dot{g}\dot{g}^{\prime }\right) }{%
\partial \theta _{h}}\right] & 0 & E^{\ast }\left[ \frac{\partial ^{2}\dot{g}%
}{\partial \theta _{h}\partial \theta ^{\prime }}\right] \\ 
0 & 0 & 0 & 0 \\ 
0 & E^{\ast }\left[ \frac{\partial ^{2}\dot{g}^{\prime }}{\partial \theta
\partial \theta _{h}}\right] & 0 & 0%
\end{array}%
\right] \left[ 
\begin{array}{c}
0 \\ 
-P\bar{g} \\ 
-P\bar{g} \\ 
-H\bar{g}%
\end{array}%
\right] =$

$\qquad =\bar{g}^{\prime }PE^{\ast }\left[ \frac{\partial \left( \dot{g}_{j}%
\dot{g}^{\prime }\right) }{\partial \theta _{h}}\right] P\bar{g}+\bar{g}%
^{\prime }PE^{\ast }\left[ \frac{\partial ^{2}\dot{g}}{\partial \theta
_{h}\partial \theta ^{\prime }}\right] H\bar{g}+\bar{g}^{\prime }H^{\prime
}E^{\ast }\left[ \frac{\partial ^{2}\dot{g}^{\prime }}{\partial \theta
\partial \theta _{h}}\right] P\bar{g}$

$\qquad =\left( \left( \bar{g}^{\prime }P\right) _{j}E^{\ast }\left[ \frac{%
\partial \left( \dot{g}_{j}\dot{g}^{\prime }\right) }{\partial \theta }%
\right] P\bar{g}\right) _{h}+\left( \left( \bar{g}^{\prime }P\right)
_{j}E^{\ast }\left[ \frac{\partial ^{2}\dot{g}_{j}}{\partial \theta \partial
\theta ^{\prime }}\right] H\bar{g}\right) _{h}+\left( \left( \bar{g}^{\prime
}H^{\prime }\right) _{j}E^{\ast }\left[ \frac{\partial ^{2}\dot{g}^{\prime }%
}{\partial \theta \partial \theta _{j}}\right] P\bar{g}\right) _{h}$

Collecting these 4 results into a single vector, we obtain:

$\bar{\Psi}_{j}\Phi _{\cdot ,jk}^{EL}\bar{\Psi}_{k}=\left[ 
\begin{array}{c}
\bar{g}^{\prime }P\bar{g} \\ 
\left( \bar{g}^{\prime }P\right) _{j}E^{\ast }\left[ \dot{g}\dot{g}_{j}\dot{g%
}^{\prime }\right] P\bar{g}+\left( \bar{g}^{\prime }P\right) _{j}E^{\ast }%
\left[ \frac{\partial \left( \dot{g}\dot{g}_{j}\right) }{\partial \theta
^{\prime }}\right] H\bar{g}+\left( \bar{g}^{\prime }H^{\prime }\right)
_{j}E^{\ast }\left[ \frac{\partial \left( \dot{g}\dot{g}^{\prime }\right) }{%
\partial \theta _{j}}\right] P\bar{g}+\left( \bar{g}^{\prime }H^{\prime
}\right) _{j}E^{\ast }\left[ \frac{\partial ^{2}\dot{g}}{\partial \theta
_{j}\partial \theta ^{\prime }}\right] H\bar{g} \\ 
\left( \bar{g}^{\prime }P\right) _{j}E^{\ast }\left[ \dot{g}_{j}\dot{g}\dot{g%
}^{\prime }\right] P\bar{g} \\ 
\left( \bar{g}^{\prime }P\right) _{j}E^{\ast }\left[ \frac{\partial \left( 
\dot{g}_{j}\dot{g}^{\prime }\right) }{\partial \theta }\right] P\bar{g}%
+\left( \bar{g}^{\prime }P\right) _{j}E^{\ast }\left[ \frac{\partial ^{2}%
\dot{g}_{j}}{\partial \theta \partial \theta ^{\prime }}\right] H\bar{g}%
+\left( \bar{g}^{\prime }H^{\prime }\right) _{j}E^{\ast }\left[ \frac{%
\partial ^{2}\dot{g}^{\prime }}{\partial \theta \partial \theta _{j}}\right]
P\bar{g}%
\end{array}%
\right] $

$\bar{\Psi}_{j}\Psi _{\cdot ,jk}^{EL}\bar{\Psi}_{k}=-\left[ 
\begin{array}{cccc}
-1 & 0 & 0 & 0 \\ 
0 & P & P & H^{\prime } \\ 
0 & P & \left( P-\Omega ^{-1}\right) & H^{\prime } \\ 
0 & H & H & -\Sigma%
\end{array}%
\right] \times $

$\qquad \times \left[ 
\begin{array}{c}
\bar{g}^{\prime }P\bar{g} \\ 
\left( \bar{g}^{\prime }P\right) _{j}E^{\ast }\left[ \dot{g}\dot{g}_{j}\dot{g%
}^{\prime }\right] P\bar{g}+\left( \bar{g}^{\prime }P\right) _{j}E^{\ast }%
\left[ \frac{\partial }{\partial \theta ^{\prime }}\left( \dot{g}\dot{g}%
_{j}\right) \right] H\bar{g}+\left( \bar{g}^{\prime }H^{\prime }\right)
_{j}E^{\ast }\left[ \frac{\partial }{\partial \theta _{j}}\left( \dot{g}\dot{%
g}^{\prime }\right) \right] P\bar{g}+\left( \bar{g}^{\prime }H^{\prime
}\right) _{j}E^{\ast }\left[ \frac{\partial ^{2}\dot{g}}{\partial \theta
_{j}\partial \theta ^{\prime }}\right] H\bar{g} \\ 
\left( \bar{g}^{\prime }P\right) _{j}E^{\ast }\left[ \dot{g}_{j}\dot{g}\dot{g%
}^{\prime }\right] P\bar{g} \\ 
\left( \bar{g}^{\prime }P\right) _{j}E^{\ast }\left[ \frac{\partial \left( 
\dot{g}_{j}\dot{g}^{\prime }\right) }{\partial \theta }\right] P\bar{g}%
+\left( \bar{g}^{\prime }P\right) _{j}E^{\ast }\left[ \frac{\partial ^{2}%
\dot{g}_{j}}{\partial \theta \partial \theta ^{\prime }}\right] H\bar{g}%
+\left( \bar{g}^{\prime }H^{\prime }\right) _{j}E^{\ast }\left[ \frac{%
\partial ^{2}\dot{g}^{\prime }}{\partial \theta \partial \theta _{j}}\right]
P\bar{g}%
\end{array}%
\right] $

$\qquad =\left[ 
\begin{array}{c}
\bar{g}^{\prime }P\bar{g} \\ 
\Xi _{1} \\ 
\Xi _{1}+\Omega ^{-1}\left( \bar{g}^{\prime }P\right) _{j}E^{\ast }\left[ 
\dot{g}_{j}\dot{g}\dot{g}^{\prime }\right] P\bar{g} \\ 
\Xi _{2}%
\end{array}%
\right] $

where

$\Xi _{1}=-P(\left( \bar{g}^{\prime }P\right) _{j}E^{\ast }\left[ \dot{g}%
\dot{g}_{j}\dot{g}^{\prime }\right] P\bar{g}+\left( \bar{g}^{\prime
}P\right) _{j}E^{\ast }\left[ \frac{\partial }{\partial \theta ^{\prime }}%
\left( \dot{g}\dot{g}_{j}\right) \right] H\bar{g}+\left( \bar{g}^{\prime
}H^{\prime }\right) _{j}E^{\ast }\left[ \frac{\partial }{\partial \theta _{j}%
}\left( \dot{g}\dot{g}^{\prime }\right) \right] P\bar{g}+$

$\qquad +\left( \bar{g}^{\prime }H^{\prime }\right) _{j}E^{\ast }\left[ 
\frac{\partial ^{2}\dot{g}}{\partial \theta _{j}\partial \theta ^{\prime }}%
\right] H\bar{g}+\left( \bar{g}^{\prime }P\right) _{j}E^{\ast }\left[ \dot{g}%
_{j}\dot{g}\dot{g}^{\prime }\right] P\bar{g})+$

$\qquad -H^{\prime }\left( \left( \bar{g}^{\prime }P\right) _{j}E^{\ast }%
\left[ \frac{\partial \left( \dot{g}_{j}\dot{g}^{\prime }\right) }{\partial
\theta }\right] P\bar{g}+\left( \bar{g}^{\prime }P\right) _{j}E^{\ast }\left[
\frac{\partial ^{2}\dot{g}_{j}}{\partial \theta \partial \theta ^{\prime }}%
\right] H\bar{g}+\left( \bar{g}^{\prime }H^{\prime }\right) _{j}E^{\ast }%
\left[ \frac{\partial ^{2}\dot{g}^{\prime }}{\partial \theta \partial \theta
_{j}}\right] P\bar{g}\right) $

$\Xi _{2}=-H(\left( \bar{g}^{\prime }P\right) _{j}E^{\ast }\left[ \dot{g}%
\dot{g}_{j}\dot{g}^{\prime }\right] P\bar{g}+\left( \bar{g}^{\prime
}P\right) _{j}E^{\ast }\left[ \frac{\partial }{\partial \theta ^{\prime }}%
\left( \dot{g}\dot{g}_{j}\right) \right] H\bar{g}+\left( \bar{g}^{\prime
}H^{\prime }\right) _{j}E^{\ast }\left[ \frac{\partial }{\partial \theta _{j}%
}\left( \dot{g}\dot{g}^{\prime }\right) \right] P\bar{g}+$

$\qquad +\left( \bar{g}^{\prime }H^{\prime }\right) _{j}E^{\ast }\left[ 
\frac{\partial ^{2}\dot{g}}{\partial \theta _{j}\partial \theta ^{\prime }}%
\right] H\bar{g}+\left( \bar{g}^{\prime }P\right) _{j}E^{\ast }\left[ \dot{g}%
_{j}\dot{g}\dot{g}^{\prime }\right] P\bar{g})+$

$\qquad +\Sigma \left( \left( \bar{g}^{\prime }P\right) _{j}E^{\ast }\left[ 
\frac{\partial \left( \dot{g}_{j}\dot{g}^{\prime }\right) }{\partial \theta }%
\right] P\bar{g}+\left( \bar{g}^{\prime }P\right) _{j}E^{\ast }\left[ \frac{%
\partial ^{2}\dot{g}_{j}}{\partial \theta \partial \theta ^{\prime }}\right]
H\bar{g}+\left( \bar{g}^{\prime }H^{\prime }\right) _{j}E^{\ast }\left[ 
\frac{\partial ^{2}\dot{g}^{\prime }}{\partial \theta \partial \theta _{j}}%
\right] P\bar{g}\right) $

Finally,

$\bar{Q}_{\cdot }^{EL}=\bar{\Psi}_{\cdot ,j}^{EL}\bar{\Psi}_{j}+\bar{\Psi}%
_{j}\Psi _{\cdot ,jk}^{EL}\bar{\Psi}_{k}/2=\left[ 
\begin{array}{c}
-\bar{g}^{\prime }P\bar{g} \\ 
P\bar{\Omega}P\bar{g}+H^{\prime }\bar{G}^{\prime }P\bar{g}+P\bar{G}H\bar{g}
\\ 
P\bar{\Omega}P\bar{g}+H^{\prime }\bar{G}^{\prime }P\bar{g}+P\bar{G}H\bar{g}
\\ 
H\bar{\Omega}P\bar{g}-\Sigma \bar{G}^{\prime }P\bar{g}+H\bar{G}H\bar{g}%
\end{array}%
\right] +\frac{1}{2}\left[ 
\begin{array}{c}
\bar{g}^{\prime }P\bar{g} \\ 
\Xi _{1} \\ 
\Xi _{1}+\Omega ^{-1}\left( \bar{g}^{\prime }P\right) _{j}E^{\ast }\left[ 
\dot{g}_{j}\dot{g}\dot{g}^{\prime }\right] P\bar{g} \\ 
\Xi _{2}%
\end{array}%
\right] $

$\qquad =\left[ 
\begin{array}{c}
-\bar{g}^{\prime }P\bar{g}/2 \\ 
\Xi _{3} \\ 
\Xi _{3}+\frac{1}{2}\Omega ^{-1}\left( \bar{g}^{\prime }P\right) _{j}E^{\ast
}\left[ \dot{g}_{j}\dot{g}\dot{g}^{\prime }\right] P\bar{g} \\ 
\Xi _{4}%
\end{array}%
\right] $

where

$\qquad \Xi _{3}=\frac{1}{2}\Xi _{1}+P\bar{\Omega}P\bar{g}+H^{\prime }\bar{G}%
^{\prime }P\bar{g}+P\bar{G}H\bar{g}$

$\qquad \Xi _{4}=\frac{1}{2}\Xi _{2}+H\bar{\Omega}P\bar{g}-\Sigma \bar{G}%
^{\prime }P\bar{g}+H\bar{G}H\bar{g}$

\section{Calculation of $\bar{R}_{\cdot }$}

We will prove that $E\left[ \left( \bar{R}_{l_{\theta }+l}^{ETEL}-\bar{R}%
_{l_{\theta }+l}^{EL}\right) \bar{\Psi}_{l_{\theta }+m}\right] =o\left(
n^{-2}\right) $. This will be done by showing that most of the term
comprised in $\bar{R}_{l_{\theta }+l}^{ETEL}-\bar{R}_{l_{\theta }+l}^{EL}$
are zero. The remaining nonzero terms will be shown to be uncorrelated with $%
\bar{\Psi}_{l_{\theta }+m}$ (up to $O\left( n^{-2}\right) $).

Since $\bar{Q}_{j}^{ETEL}=\bar{Q}_{j}^{EL}\equiv \bar{Q}_{j}$%
\begin{eqnarray*}
R_{l}^{ETEL}-R_{l}^{EL} &=&\left( \bar{\Psi}_{l,j}^{ETEL}-\bar{\Psi}%
_{l,j}^{EL}\right) \bar{Q}_{j}+\left( \Psi _{l,jk}^{ETEL}-\Psi
_{l,jk}^{EL}\right) \bar{\Psi}_{k}\bar{Q}_{j}+ \\
&&+\frac{1}{2}\left( \bar{\Psi}_{l,jk}^{ETEL}-\bar{\Psi}_{l,jk}^{EL}\right) 
\bar{\Psi}_{j}\bar{\Psi}_{k}+\frac{1}{6}\left( \Psi _{l,jkh}^{ETEL}-\Psi
_{l,jkh}^{EL}\right) \bar{\Psi}_{j}\bar{\Psi}_{k}\bar{\Psi}_{h}
\end{eqnarray*}

\subsection{Calculation of $\left( \bar{\Psi}_{l_{\protect\theta %
}+l,j}^{ETEL}-\bar{\Psi}_{l_{\protect\theta }+l,j}^{EL}\right) \bar{Q}_{j}$%
\label{Rterm1}}

$\left( \bar{\Psi}_{l_{\theta }+l,\cdot }^{ETEL}-\bar{\Psi}_{l_{\theta
}+l,\cdot }^{EL}\right) =\Phi _{l_{\theta }+l,h}\left( \bar{\Phi}_{h,\cdot
}^{ETEL}-\bar{\Phi}_{h,\cdot }^{EL}\right) =$

$\qquad =-\left[ 
\begin{array}{cccc}
0 & H_{l\cdot } & H_{l\cdot } & -\Sigma _{l\cdot }%
\end{array}%
\right] \left( \left[ 
\begin{array}{cccc}
-1 & 0 & \bar{g}^{\prime } & 0 \\ 
0 & 0 & \bar{g}\bar{g}^{\prime } & \bar{G} \\ 
\bar{g} & \bar{g}\bar{g}^{\prime } & -\bar{g}\bar{g}^{\prime } & 0 \\ 
0 & \bar{G}^{\prime } & 0 & 0%
\end{array}%
\right] -\left[ 
\begin{array}{cccc}
-1 & 0 & \bar{g}^{\prime } & 0 \\ 
0 & 0 & \bar{g}\bar{g}^{\prime } & \bar{G} \\ 
0 & \bar{g}\bar{g}^{\prime } & -\bar{g}\bar{g}^{\prime } & 0 \\ 
0 & \bar{G}^{\prime } & 0 & 0%
\end{array}%
\right] \right) =$

$\qquad =\left[ 
\begin{array}{cccc}
-H_{l\cdot }\bar{g} & 0 & 0 & 0%
\end{array}%
\right] $

$\left( \bar{\Psi}_{l_{\theta }+l,j}^{ETEL}-\bar{\Psi}_{l_{\theta
}+l,j}^{EL}\right) \bar{Q}_{j}=\left( \bar{\Psi}_{l_{\theta }+l,1}^{ETEL}-%
\bar{\Psi}_{l_{\theta }+l,1}^{EL}\right) \bar{Q}_{1}=-H_{l\cdot }\bar{g}%
\left( -\bar{g}^{\prime }P\bar{g}/2\right) =H_{l\cdot }\bar{g}\bar{g}%
^{\prime }P\bar{g}/2$

\subsection{Calculation of $\left( \Psi _{l,jk}^{ETEL}-\Psi
_{l,jk}^{EL}\right) \bar{\Psi}_{k}\bar{Q}_{j}$\label{Rterm2}}

$\left( \Psi _{l,jk}^{ETEL}-\Psi _{l,jk}^{EL}\right) \bar{\Psi}_{k}\bar{Q}%
_{j}=\Phi _{lh}^{-1}\left( \Phi _{h,jk}^{ETEL}-\Phi _{h,jk}^{EL}\right) \bar{%
\Psi}_{k}\bar{Q}_{j}$ where $\left( \Phi _{h,jk}^{ETEL}-\Phi
_{h,jk}^{EL}\right) \bar{\Psi}_{k}$ can be calculated as in Section \ref%
{secpsiquad}:

Since $\left( \dot{\phi}^{ETEL}-\dot{\phi}^{EL}\right) =\left[ 
\begin{array}{c}
0 \\ 
0 \\ 
\left( \tau -1\right) \dot{g}+\left( \dot{\tau}-\dot{\varepsilon}\right) 
\dot{g}\dot{g}^{\prime }\kappa \\ 
\left( \dot{\tau}-\dot{\varepsilon}\right) \dot{G}^{\prime }\kappa +\dot{\tau%
}\dot{G}^{\prime }\lambda \dot{g}^{\prime }\kappa -\left( \dot{\tau}-\tau
\right) \dot{G}^{\prime }\lambda%
\end{array}%
\right] $

$\left( \Phi _{l_{\tau }+h,\cdot k}^{ETEL}-\Phi _{l_{\tau }+h,\cdot
k}^{EL}\right) \bar{\Psi}_{k}=\mathbf{0}$

$\left( \Phi _{l_{\kappa }+h,\cdot k}^{ETEL}-\Phi _{l_{\kappa }+h,\cdot
k}^{EL}\right) \bar{\Psi}_{k}=\mathbf{0}$

$\left( \Phi _{l_{\lambda }+h,\cdot k}^{ETEL}-\Phi _{l_{\lambda }+h,\cdot
k}^{EL}\right) \bar{\Psi}_{k}=\left[ 
\begin{array}{cccc}
0 & 0 & 0 & G_{h\cdot } \\ 
0 & -2E^{\ast }\left[ \dot{g}\dot{g}^{\prime }\dot{g}_{h}\right] & E^{\ast } 
\left[ \dot{g}_{h}\dot{g}\dot{g}^{\prime }\right] & 0 \\ 
0 & E^{\ast }\left[ \dot{g}_{h}\dot{g}\dot{g}^{\prime }\right] & 0 & 0 \\ 
G_{h\cdot }^{\prime } & 0 & 0 & 0%
\end{array}%
\right] \left[ 
\begin{array}{c}
0 \\ 
-P\bar{g} \\ 
-P\bar{g} \\ 
-H\bar{g}%
\end{array}%
\right] \allowbreak =\left[ 
\begin{array}{c}
-G_{h\cdot }H\bar{g} \\ 
E^{\ast }\left[ \dot{g}\dot{g}^{\prime }\dot{g}_{h}\right] P\bar{g} \\ 
-E^{\ast }\left[ \dot{g}\dot{g}^{\prime }\dot{g}_{h}\right] P\bar{g} \\ 
0%
\end{array}%
\right] $

$\left( \Phi _{l_{\theta }+h,\cdot k}^{ETEL}-\Phi _{l_{\theta }+h,\cdot
k}^{EL}\right) \bar{\Psi}_{k}=\left[ 
\begin{array}{cccc}
0 & 0 & G_{\cdot h}^{\prime } & 0 \\ 
0 & -E^{\ast }\left[ \frac{\partial \left( \dot{g}\dot{g}^{\prime }\right) }{%
\partial \theta _{h}}\right] & E^{\ast }\left[ \frac{\partial \left( \dot{g}%
\dot{g}^{\prime }\right) }{\partial \theta _{h}}\right] & 0 \\ 
G_{\cdot h} & E^{\ast }\left[ \frac{\partial \left( \dot{g}\dot{g}^{\prime
}\right) }{\partial \theta _{h}}\right] & -E^{\ast }\left[ \frac{\partial
\left( \dot{g}\dot{g}^{\prime }\right) }{\partial \theta _{h}}\right] & 0 \\ 
0 & 0 & 0 & 0%
\end{array}%
\right] \left[ 
\begin{array}{c}
0 \\ 
-P\bar{g} \\ 
-P\bar{g} \\ 
-H\bar{g}%
\end{array}%
\right] \allowbreak =\left[ 
\begin{array}{c}
-G_{\cdot h}^{\prime }P\bar{g} \\ 
0 \\ 
0 \\ 
0%
\end{array}%
\right] =\left[ 
\begin{array}{c}
0 \\ 
0 \\ 
0 \\ 
0%
\end{array}%
\right] $

Then, since $\left( \Phi _{h,jk}^{ETEL}-\Phi _{h,jk}^{EL}\right) \bar{\Psi}%
_{k}\bar{Q}_{j}=\bar{Q}_{\cdot }^{\prime }\left( \Phi _{h,\cdot
k}^{ETEL}-\Phi _{h,\cdot k}^{EL}\right) \bar{\Psi}_{k}$, we obtain

$\bar{Q}_{\cdot }^{\prime }\left( \Phi _{l_{\tau }+h,\cdot k}^{ETEL}-\Phi
_{l_{\tau }+h,\cdot k}^{EL}\right) \bar{\Psi}_{k}=0$

$\bar{Q}_{\cdot }^{\prime }\left( \Phi _{l_{\kappa }+h,\cdot k}^{ETEL}-\Phi
_{l_{\kappa }+h,\cdot k}^{EL}\right) \bar{\Psi}_{k}=0$

$\bar{Q}_{\cdot }^{\prime }\left( \Phi _{l_{\lambda }+h,\cdot k}^{ETEL}-\Phi
_{l_{\lambda }+h,\cdot k}^{EL}\right) \bar{\Psi}_{k}=\left[ 
\begin{array}{cccc}
-\bar{g}^{\prime }P\bar{g}/2 & \Xi _{3}^{\prime } & \Xi _{3}^{\prime }+\frac{%
1}{2}\left( \bar{g}^{\prime }P\right) _{j}\bar{g}^{\prime }PE^{\ast }\left[ 
\dot{g}_{j}\dot{g}\dot{g}^{\prime }\right] \Omega ^{-1} & \Xi _{4}^{\prime }%
\end{array}%
\right] \left[ 
\begin{array}{c}
-G_{h\cdot }H\bar{g} \\ 
E^{\ast }\left[ \dot{g}\dot{g}^{\prime }\dot{g}_{h}\right] P\bar{g} \\ 
-E^{\ast }\left[ \dot{g}\dot{g}^{\prime }\dot{g}_{h}\right] P\bar{g} \\ 
0%
\end{array}%
\right] $

$\qquad =\frac{1}{2}\bar{g}^{\prime }P\bar{g}G_{h\cdot }H\bar{g}-\frac{1}{2}%
\left( \bar{g}^{\prime }P\right) _{j}\bar{g}^{\prime }PE^{\ast }\left[ \dot{g%
}_{j}\dot{g}\dot{g}^{\prime }\right] \Omega ^{-1}E^{\ast }\left[ \dot{g}\dot{%
g}^{\prime }\dot{g}_{h}\right] P\bar{g}=$

$\qquad =\frac{1}{2}G_{h\cdot }H\bar{g}\bar{g}^{\prime }P\bar{g}-\frac{1}{2}%
E^{\ast }\left[ \dot{g}_{h}\dot{g}_{k}\dot{g}^{\prime }\right] \Omega
^{-1}E^{\ast }\left[ \dot{g}_{j}\dot{g}\dot{g}^{\prime }\right] P\bar{g}%
\left( \bar{g}^{\prime }P\right) _{j}\left( \bar{g}^{\prime }P\right) _{k}$

$\bar{Q}_{\cdot }^{\prime }\left( \Phi _{l_{\theta }+h,\cdot k}^{ETEL}-\Phi
_{l_{\theta }+h,\cdot k}^{EL}\right) \bar{\Psi}_{k}=0$

Collecting the previous four quantities into a single vector, we have

$\left( \Phi _{\cdot ,jk}^{ETEL}-\Phi _{\cdot ,\cdot jk}^{EL}\right) \bar{Q}%
_{j}\bar{\Psi}_{k}=\left[ 
\begin{array}{c}
0 \\ 
0 \\ 
\frac{1}{2}GH\bar{g}\bar{g}^{\prime }P\bar{g}+E^{\ast }\left[ \dot{g}\dot{g}%
_{k}\dot{g}^{\prime }\right] \Omega ^{-1}E^{\ast }\left[ \dot{g}_{j}\dot{g}%
\dot{g}^{\prime }\right] P\bar{g}\left( \bar{g}^{\prime }P\right) _{j}\left( 
\bar{g}^{\prime }P\right) _{k} \\ 
0%
\end{array}%
\right] $

$\Phi _{l_{\theta }+l,\cdot }^{-1}\left( \Phi _{\cdot ,jk}^{ETEL}-\Phi
_{\cdot ,\cdot jk}^{EL}\right) \bar{Q}_{j}\bar{\Psi}_{k}=$

$=-\left[ 
\begin{array}{cccc}
0 & H_{l\cdot } & H_{l\cdot } & -\Sigma _{l\cdot }%
\end{array}%
\right] \left[ 
\begin{array}{c}
0 \\ 
0 \\ 
\frac{1}{2}GH\bar{g}\bar{g}^{\prime }P\bar{g}+E^{\ast }\left[ \dot{g}\dot{g}%
_{k}\dot{g}^{\prime }\right] \Omega ^{-1}E^{\ast }\left[ \dot{g}_{j}\dot{g}%
\dot{g}^{\prime }\right] P\bar{g}\left( \bar{g}^{\prime }P\right) _{j}\left( 
\bar{g}^{\prime }P\right) _{k} \\ 
0%
\end{array}%
\right] $

$=-\left( \frac{1}{2}HGH\bar{g}\bar{g}^{\prime }P\bar{g}\right) _{l}-\left(
HE^{\ast }\left[ \dot{g}\dot{g}_{k}\dot{g}^{\prime }\right] \Omega
^{-1}E^{\ast }\left[ \dot{g}_{j}\dot{g}\dot{g}^{\prime }\right] P\bar{g}%
\left( \bar{g}^{\prime }P\right) _{j}\left( \bar{g}^{\prime }P\right)
_{k}\right) _{l}$

$=-\frac{1}{2}\left( \left( \left( G^{\prime }\Omega ^{-1}G\right)
^{-1}G^{\prime }\Omega ^{-1}\right) GH\bar{g}\bar{g}^{\prime }P\bar{g}%
\right) _{l}-\left( HE^{\ast }\left[ \dot{g}\dot{g}_{k}\dot{g}^{\prime }%
\right] \Omega ^{-1}E^{\ast }\left[ \dot{g}_{j}\dot{g}\dot{g}^{\prime }%
\right] P\bar{g}\left( \bar{g}^{\prime }P\right) _{j}\left( \bar{g}^{\prime
}P\right) _{k}\right) _{l}$

$=-\left( \frac{1}{2}H\bar{g}\bar{g}^{\prime }P\bar{g}\right) _{l}-\left(
HE^{\ast }\left[ \dot{g}\dot{g}_{k}\dot{g}^{\prime }\right] \Omega
^{-1}E^{\ast }\left[ \dot{g}_{j}\dot{g}\dot{g}_{m}\right] \left( P\bar{g}%
\right) _{m}\left( P\bar{g}\right) _{j}\left( P\bar{g}\right) _{k}\right)
_{l}$

The first term is exactly cancels the first term of the expansion of $\bar{R}%
_{\cdot }$ (see section \ref{Rterm1}).

The second term is uncorrelated with $\bar{\Psi}_{l_{\theta }+l}$, which
implies that it does not contribute to the $n^{-2}$ variance. Indeed, all
the random contributions to the second term are of the form $P\bar{g}$,
while $\bar{\Psi}_{l_{\theta }+l}$ is of the form $H\bar{g}$ and $E^{\ast }%
\left[ P\bar{g}\bar{g}^{\prime }H\right] =P\Omega H=0$.

\subsection{Calculation of $\bar{\Psi}_{j}\left( \bar{\Psi}_{l,jk}^{ETEL}-%
\bar{\Psi}_{l,jk}^{EL}\right) \bar{\Psi}_{k}$\label{Rterm3}}

$\bar{\Psi}_{j}\left( \bar{\Psi}_{l,jk}^{ETEL}-\bar{\Psi}_{l,jk}^{EL}\right) 
\bar{\Psi}_{k}=\Phi _{lh}^{-1}\bar{\Psi}_{j}\left( \bar{\Phi}_{h,jk}^{ETEL}-%
\bar{\Phi}_{h,jk}^{EL}\right) \bar{\Psi}_{k}$

Since $\left( \dot{\phi}^{ETEL}-\dot{\phi}^{EL}\right) =\left[ 
\begin{array}{c}
0 \\ 
0 \\ 
\left( \tau -1\right) \dot{g}+\left( \dot{\tau}-\dot{\varepsilon}\right) 
\dot{g}\dot{g}^{\prime }\kappa \\ 
\left( \dot{\tau}-\dot{\varepsilon}\right) \dot{G}^{\prime }\kappa +\dot{\tau%
}\dot{G}^{\prime }\lambda \dot{g}^{\prime }\kappa -\left( \dot{\tau}-\tau
\right) \dot{G}^{\prime }\lambda%
\end{array}%
\right] ,$

we have, in analogy with Section \ref{secpsiquad},

$\bar{\Psi}_{j}\left( \bar{\Phi}_{l_{\tau }+h,jk}^{ETEL}-\bar{\Phi}_{l_{\tau
}+h,jk}^{EL}\right) \bar{\Psi}_{k}=\mathbf{0}$

$\bar{\Psi}_{j}\left( \bar{\Phi}_{l_{\kappa }+h,jk}^{ETEL}-\bar{\Phi}%
_{l_{\kappa }+h,jk}^{EL}\right) \bar{\Psi}_{k}=\mathbf{0}$

$\bar{\Psi}_{j}\left( \bar{\Phi}_{l_{\lambda }+h,jk}^{ETEL}-\bar{\Phi}%
_{l_{\lambda }+h,jk}^{EL}\right) \bar{\Psi}_{k}=\left[ 
\begin{array}{cccc}
0 & -\bar{g}^{\prime }P & -\bar{g}^{\prime }P & -\bar{g}^{\prime }H^{\prime }%
\end{array}%
\right] \left[ 
\begin{array}{cccc}
0 & 0 & 0 & \bar{G}_{h\cdot } \\ 
0 & -2\overline{g_{h}gg^{\prime }} & \overline{g_{h}gg^{\prime }} & 0 \\ 
0 & \overline{g_{h}gg^{\prime }} & 0 & 0 \\ 
\bar{G}_{h\cdot }^{\prime } & 0 & 0 & 0%
\end{array}%
\right] \left[ 
\begin{array}{c}
0 \\ 
-P\bar{g} \\ 
-P\bar{g} \\ 
-H\bar{g}%
\end{array}%
\right] $

$\qquad =\mathbf{0}$

$\bar{\Psi}_{j}\left( \bar{\Phi}_{l_{\theta }+h,jk}^{ETEL}-\bar{\Phi}%
_{l_{\theta }+h,jk}^{EL}\right) \bar{\Psi}_{k}=\left[ 
\begin{array}{cccc}
0 & -\bar{g}^{\prime }P & -\bar{g}^{\prime }P & -\bar{g}^{\prime }H^{\prime }%
\end{array}%
\right] \left[ 
\begin{array}{cccc}
0 & 0 & \bar{G}_{\cdot h}^{\prime } & 0 \\ 
0 & -\left( \overline{\frac{\partial \left( gg^{\prime }\right) }{\partial
\theta _{h}}}\right) & \left( \overline{\frac{\partial \left( gg^{\prime
}\right) }{\partial \theta _{h}}}\right) & 0 \\ 
\bar{G}_{\cdot h} & \left( \overline{\frac{\partial \left( gg^{\prime
}\right) }{\partial \theta _{h}}}\right) & -\left( \overline{\frac{\partial
\left( gg^{\prime }\right) }{\partial \theta _{h}}}\right) & 0 \\ 
0 & 0 & 0 & 0%
\end{array}%
\right] \left[ 
\begin{array}{c}
0 \\ 
-P\bar{g} \\ 
-P\bar{g} \\ 
-H\bar{g}%
\end{array}%
\right] $

$\qquad =\mathbf{0}$

where $\overline{g_{h}gg^{\prime }}\equiv \left[ n^{-1/2}\sum_{i}\left( \dot{%
g}_{h}\dot{g}\dot{g}^{\prime }-E\left[ \dot{g}_{h}\dot{g}\dot{g}^{\prime }%
\right] \right) \right] _{\beta =\beta ^{\ast }}$ and $\overline{\frac{%
\partial \left( gg^{\prime }\right) }{\partial \theta _{h}}}=\left[
n^{-1/2}\sum_{i}\left( \overline{\frac{\partial \left( gg^{\prime }\right) }{%
\partial \theta _{h}}}-E\left[ \frac{\partial \left( gg^{\prime }\right) }{%
\partial \theta _{h}}\right] \right) \right] _{\beta =\beta ^{\ast }}.$

It follows that $\bar{\Psi}_{j}\left( \bar{\Psi}_{l,jk}^{ETEL}-\bar{\Psi}%
_{l,jk}^{EL}\right) \bar{\Psi}_{k}=\Phi _{lh}^{-1}\bar{\Psi}_{j}\left( \bar{%
\Phi}_{h,jk}^{ETEL}-\bar{\Phi}_{h,jk}^{EL}\right) \bar{\Psi}_{k}=0$.

\subsection{Calculation of $\left( \Psi _{l,jkh}^{ETEL}-\Psi
_{l,jkh}^{EL}\right) \bar{\Psi}_{j}\bar{\Psi}_{k}\bar{\Psi}_{h}$\label%
{Rterm4}}

\subsubsection{Simplification}

Since $\bar{R}_{l_{\theta }+l}$ contributes to the $n^{-2}$ variance only
through its correlation with $\bar{\Psi}_{l_{\theta }+m}$, we can omit terms
of $\bar{R}_{l_{\theta }+l}\,\ $that are uncorrelated with $\bar{\Psi}%
_{l_{\theta }+m}$. The $O\left( n^{-2}\right) $ correlation of $\left( \Psi
_{l_{\theta }+l,jkq}^{ETEL}-\Psi _{l_{\theta }+l,jkq}^{EL}\right) \bar{\Psi}%
_{j}\bar{\Psi}_{k}\bar{\Psi}_{q}$ with $\bar{\Psi}_{l_{\theta }+m}$ can be
written as

$\left( \Psi _{l_{\theta }+l,jkq}^{ETEL}-\Psi _{l_{\theta
}+l,jkq}^{EL}\right) E\left[ \bar{\Psi}_{j}\bar{\Psi}_{k}\right] E\left[ 
\bar{\Psi}_{q}\bar{\Psi}_{l_{\theta }+m}\right] +\left( \Psi _{l_{\theta
}+l,jkq}^{ETEL}-\Psi _{l_{\theta }+l,jkq}^{EL}\right) E\left[ \bar{\Psi}_{j}%
\bar{\Psi}_{q}\right] E\left[ \bar{\Psi}_{k}\bar{\Psi}_{l_{\theta }+m}\right]
+$

$\qquad +\left( \Psi _{l_{\theta }+l,jkq}^{ETEL}-\Psi _{l_{\theta
}+l,jkq}^{EL}\right) E\left[ \bar{\Psi}_{j}\bar{\Psi}_{l_{\theta }+m}\right]
E\left[ \bar{\Psi}_{k}\bar{\Psi}_{q}\right] $.

Since

$E\left[ \bar{\Psi}_{\cdot }\bar{\Psi}_{\cdot }^{\prime }\right] =\left[ 
\begin{array}{cccc}
0 & 0 & 0 & 0 \\ 
0 & P & P & 0 \\ 
0 & P & P & 0 \\ 
0 & 0 & 0 & \Sigma%
\end{array}%
\right] ,$

a term of the form $E\left[ \bar{\Psi}_{q}\bar{\Psi}_{l_{\theta }+m}\right] $
will be nonzero only if $q>l_{\theta }$. It follows that, in the sum $\left(
\Psi _{l,jkq}^{ETEL}-\Psi _{l,jkq}^{EL}\right) \bar{\Psi}_{j}\bar{\Psi}_{k}%
\bar{\Psi}_{q}$, we only need to keep terms such that either $j>l_{\theta }$
or $k>l_{\theta }$ or $q>l_{\theta }$.

Letting $\Psi _{l_{\theta }+l,jkq}^{ETEL-EL}=\left( \Psi _{l_{\theta
}+l,jkq}^{ETEL}-\Psi _{l_{\theta }+l,jkq}^{EL}\right) $, and making use of
the fact that $\Psi _{l,jkq}=\Psi _{l,kjq}=\Psi _{l,qjk}=\ldots $, we have,

\begin{eqnarray*}
&&\sum_{j,k,q}\Psi _{l_{\theta }+l,jkq}^{ETEL-EL}\bar{\Psi}_{j}\bar{\Psi}_{k}%
\bar{\Psi}_{q}1\left( j>l_{\theta }\text{ or }k>l_{\theta }\text{ or }%
q>l_{\theta }\right) \\
&=&\sum_{j}\sum_{k}\sum_{q>l_{\theta }}\Psi _{l_{\theta }+l,jkq}^{ETEL-EL}%
\bar{\Psi}_{j}\bar{\Psi}_{k}\bar{\Psi}_{q}+\sum_{j}\sum_{k}\sum_{q\leq
l_{\theta }}\Psi _{l_{\theta }+l,jkq}^{ETEL-EL}\bar{\Psi}_{j}\bar{\Psi}_{k}%
\bar{\Psi}_{q}1\left( j>l_{\theta }\text{ or }k>l_{\theta }\right) \\
&=&\sum_{j}\sum_{k}\sum_{q>l_{\theta }}\Psi _{l_{\theta }+l,jkq}^{ETEL-EL}%
\bar{\Psi}_{j}\bar{\Psi}_{k}\bar{\Psi}_{q}+\sum_{j}\sum_{k\leq l_{\theta
}}\sum_{q}\Psi _{l_{\theta }+l,jkq}^{ETEL-EL}\bar{\Psi}_{j}\bar{\Psi}_{k}%
\bar{\Psi}_{q}1\left( j>l_{\theta }\text{ or }q>l_{\theta }\right) \\
&=&\sum_{j}\sum_{k}\sum_{q>l_{\theta }}\Psi _{l_{\theta }+l,jkq}^{ETEL-EL}%
\bar{\Psi}_{j}\bar{\Psi}_{k}\bar{\Psi}_{q}+\sum_{j}\sum_{k\leq l_{\theta
}}\sum_{q>l_{\theta }}\Psi _{l_{\theta }+l,jkq}^{ETEL-EL}\bar{\Psi}_{j}\bar{%
\Psi}_{k}\bar{\Psi}_{q}+ \\
&&+\sum_{j}\sum_{k\leq l_{\theta }}\sum_{q\leq l_{\theta }}\Psi _{l_{\theta
}+l,jkq}^{ETEL-EL}\bar{\Psi}_{j}\bar{\Psi}_{k}\bar{\Psi}_{q}1\left(
j>l_{\theta }\right) \\
&=&\sum_{j}\sum_{k}\sum_{q>l_{\theta }}\Psi _{l_{\theta }+l,jkq}^{ETEL-EL}%
\bar{\Psi}_{j}\bar{\Psi}_{k}\bar{\Psi}_{q}+\sum_{j}\sum_{k\leq l_{\theta
}}\sum_{q>l_{\theta }}\Psi _{l_{\theta }+l,jkq}^{ETEL-EL}\bar{\Psi}_{j}\bar{%
\Psi}_{k}\bar{\Psi}_{q}+\sum_{j>l_{\theta }}\sum_{k\leq l_{\theta
}}\sum_{q\leq l_{\theta }}\Psi _{l_{\theta }+l,jkq}^{ETEL-EL}\bar{\Psi}_{j}%
\bar{\Psi}_{k}\bar{\Psi}_{q} \\
&=&\sum_{j}\sum_{k}\sum_{q>l_{\theta }}\Psi _{l_{\theta }+l,jkq}^{ETEL-EL}%
\bar{\Psi}_{j}\bar{\Psi}_{k}\bar{\Psi}_{q}+\sum_{j}\sum_{k\leq l_{\theta
}}\sum_{q>l_{\theta }}\Psi _{l_{\theta }+l,jkq}^{ETEL-EL}\bar{\Psi}_{j}\bar{%
\Psi}_{k}\bar{\Psi}_{q}+\sum_{j\leq l_{\theta }}\sum_{k\leq l_{\theta
}}\sum_{q>l_{\theta }}\Psi _{l_{\theta }+l,jkq}^{ETEL-EL}\bar{\Psi}_{j}\bar{%
\Psi}_{k}\bar{\Psi}_{q} \\
&=&\sum_{j}\sum_{k}\sum_{q>l_{\theta }}\Psi _{l_{\theta }+l,jkq}^{ETEL-EL}%
\bar{\Psi}_{j}\bar{\Psi}_{k}\bar{\Psi}_{q}+\frac{1}{2}\sum_{j\leq l_{\theta
}}\sum_{k}\sum_{q>l_{\theta }}\Psi _{l_{\theta }+l,jkq}^{ETEL-EL}\bar{\Psi}%
_{j}\bar{\Psi}_{k}\bar{\Psi}_{q}+\frac{1}{2}\sum_{j}\sum_{k\leq l_{\theta
}}\sum_{q>l_{\theta }}\Psi _{l_{\theta }+l,jkq}^{ETEL-EL}\bar{\Psi}_{j}\bar{%
\Psi}_{k}\bar{\Psi}_{q}+ \\
&&+\sum_{j\leq l_{\theta }}\sum_{k\leq l_{\theta }}\sum_{q>l_{\theta }}\Psi
_{l_{\theta }+l,jkq}^{ETEL-EL}\bar{\Psi}_{j}\bar{\Psi}_{k}\bar{\Psi}_{q} \\
&=&\sum_{j}\sum_{k}\sum_{q>l_{\theta }}\Psi _{l_{\theta }+l,jkq}^{ETEL-EL}%
\bar{\Psi}_{j}\bar{\Psi}_{k}\bar{\Psi}_{q}\xi _{jk}
\end{eqnarray*}%
where%
\[
\xi _{jk}=\left\{ 
\begin{array}{ll}
1 & \text{if }j>l_{\theta }\text{ and }k>l_{\theta } \\ 
3/2 & \text{if }\left( j>l_{\theta }\text{ and }k\leq l_{\theta }\right) 
\text{ or }\left( j\leq l_{\theta }\text{ and }k>l_{\theta }\right) \\ 
3 & \text{if }j\leq l_{\theta }\text{ and }k\leq l_{\theta }%
\end{array}%
\right.
\]%
$\sum_{j}\sum_{k}\sum_{q>l_{\theta }}\Psi _{l_{\theta }+l,jkq}^{ETEL-EL}\bar{%
\Psi}_{j}\bar{\Psi}_{k}\bar{\Psi}_{q}\xi _{jk}=\Phi _{l_{\theta
}+l,h}^{-1}\sum_{q>l_{\theta }}\left( \sum_{j}\sum_{k}\left( \Phi
_{h,jkq}^{ETEL}-\Phi _{h,jkq}^{EL}\right) \xi _{jk}\bar{\Psi}_{j}\bar{\Psi}%
_{k}\right) \bar{\Psi}_{q}$

\subsubsection{Calculation of $\left( \Phi _{h,~j~k~l_{\protect\theta %
}+q}^{ETEL}-\Phi _{h,~j~k~l_{\protect\theta }+q}^{EL}\right) \protect\xi %
_{jk}\bar{\Psi}_{j}\bar{\Psi}_{k}$}

Noting that $\left( \dot{\phi}^{ETEL}-\dot{\phi}^{EL}\right) =\left[ 
\begin{array}{c}
0 \\ 
0 \\ 
\left( \tau -1\right) \dot{g}+\left( \dot{\tau}-\dot{\varepsilon}\right) 
\dot{g}\dot{g}^{\prime }\kappa \\ 
\left( \dot{\tau}-\dot{\varepsilon}\right) \dot{G}^{\prime }\kappa +\dot{\tau%
}\dot{G}^{\prime }\lambda \dot{g}^{\prime }\kappa -\left( \dot{\tau}-\tau
\right) \dot{G}^{\prime }\lambda%
\end{array}%
\right] $, we obtain:

(When writing the intermediate steps of the calculations, we omit the terms
that will be multiplied by $\lambda $ or $\kappa $ after all the derivatives
have been evaluated, since these terms will vanish when the true values $%
\lambda =0$ and $\kappa =0$ are substituted in.)

$\left( \Phi _{l_{\tau }+h,jk\left( l_{\theta }+q\right) }^{ETEL}-\Phi
_{l_{\tau }+h,jk\left( l_{\theta }+q\right) }^{EL}\right) =\mathbf{0}$

$\left( \Phi _{l_{\kappa }+h,jk\left( l_{\theta }+q\right) }^{ETEL}-\Phi
_{l_{\kappa }+h,jk\left( l_{\theta }+q\right) }^{EL}\right) =\mathbf{0}$

$\left( \Phi _{l_{\lambda }+h,jk\left( l_{\theta }+q\right) }^{ETEL}-\Phi
_{l_{\lambda }+h,jk\left( l_{\theta }+q\right) }^{EL}\right) =E^{\ast }\left[
\frac{\partial ^{3}}{\partial \beta \partial \beta ^{\prime }\partial \theta
_{q}}\left( \left( \tau -1\right) \dot{g}_{h}\right) +\frac{\partial ^{3}}{%
\partial \beta \partial \beta ^{\prime }\partial \theta _{q}}\left( \left( 
\dot{\tau}-\dot{\varepsilon}\right) \dot{g}_{h}\dot{g}^{\prime }\kappa
\right) \right] =$

$\qquad =\left[ 
\begin{array}{cccc}
0 & 0 & 0 & E^{\ast }\left[ \frac{\partial ^{2}\dot{g}_{h}}{\partial \theta
^{\prime }\partial \theta _{q}}\right] \\ 
0 & 0 & 0 & 0 \\ 
0 & 0 & 0 & 0 \\ 
E^{\ast }\left[ \frac{\partial ^{2}\dot{g}_{h}^{\prime }}{\partial \theta
\partial \theta _{q}}\right] & 0 & 0 & 0%
\end{array}%
\right] +$

$\qquad +\left[ 
\begin{array}{cccc}
0 & 0 & 0 & 0 \\ 
0 & -2E^{\ast }\left[ \frac{\partial }{\partial \theta _{q}}\left( \dot{%
\varepsilon}^{2}\dot{g}_{h}\dot{g}\dot{g}^{\prime }\right) \right] & E^{\ast
}\left[ \frac{\partial }{\partial \theta _{q}}\left( \dot{\tau}\dot{g}_{h}%
\dot{g}\dot{g}^{\prime }\right) \right] & 0 \\ 
0 & E^{\ast }\left[ \frac{\partial }{\partial \theta _{q}}\left( \dot{\tau}%
\dot{g}_{h}\dot{g}\dot{g}^{\prime }\right) \right] & 0 & 0 \\ 
0 & 0 & 0 & 0%
\end{array}%
\right] $

\qquad $=\left[ 
\begin{array}{cccc}
0 & 0 & 0 & E^{\ast }\left[ \frac{\partial ^{2}\dot{g}_{h}}{\partial \theta
^{\prime }\partial \theta _{q}}\right] \\ 
0 & -2E^{\ast }\left[ \frac{\partial \left( \dot{g}_{h}\dot{g}\dot{g}%
^{\prime }\right) }{\partial \theta _{q}}\right] & E^{\ast }\left[ \frac{%
\partial \left( \dot{g}_{h}\dot{g}\dot{g}^{\prime }\right) }{\partial \theta
_{q}}\right] & 0 \\ 
0 & E^{\ast }\left[ \frac{\partial \left( \dot{g}_{h}\dot{g}\dot{g}^{\prime
}\right) }{\partial \theta _{q}}\right] & 0 & 0 \\ 
E^{\ast }\left[ \frac{\partial ^{2}\dot{g}_{h}^{\prime }}{\partial \theta
\partial \theta _{q}}\right] & 0 & 0 & 0%
\end{array}%
\right] $

$\left( \Phi _{l_{\lambda }+h,jk\left( l_{\theta }+q\right) }^{ETEL}-\Phi
_{l_{\lambda }+h,jk\left( l_{\theta }+q\right) }^{EL}\right) =$

$=E^{\ast }\left[ \frac{\partial ^{3}}{\partial \beta \partial \beta
^{\prime }\partial \theta _{q}}\left( \left( \dot{\tau}-\dot{\varepsilon}%
\right) \dot{G}_{\cdot h}^{\prime }\kappa \right) +\frac{\partial ^{3}}{%
\partial \beta \partial \beta ^{\prime }\partial \theta _{q}}\left( \dot{\tau%
}\dot{G}_{\cdot h}^{\prime }\lambda \dot{g}^{\prime }\kappa \right) +\frac{%
\partial ^{3}}{\partial \beta \partial \beta ^{\prime }\partial \theta _{q}}%
\left( \left( \tau -\dot{\tau}\right) \dot{G}_{\cdot h}^{\prime }\lambda
\right) \right] =$

$\qquad =\left[ 
\begin{array}{cccc}
0 & 0 & 0 & 0 \\ 
0 & -E^{\ast }\left[ \left( \frac{\partial }{\partial \theta _{q}}\left( 
\dot{\varepsilon}^{2}\dot{g}\dot{G}_{\cdot h}^{\prime }+\dot{\varepsilon}^{2}%
\dot{G}_{\cdot h}\dot{g}^{\prime }\right) \right) \right] & E^{\ast }\left[ 
\frac{\partial }{\partial \theta _{q}}\left( \dot{\tau}\dot{G}_{\cdot h}\dot{%
g}^{\prime }\right) \right] & 0 \\ 
0 & E^{\ast }\left[ \frac{\partial }{\partial \theta _{q}}\left( \dot{\tau}%
\dot{g}\dot{G}_{\cdot h}^{\prime }\right) \right] & 0 & 0 \\ 
0 & 0 & 0 & 0%
\end{array}%
\right] +$

$\qquad +\left[ 
\begin{array}{cccc}
0 & 0 & 0 & 0 \\ 
0 & 0 & E^{\ast }\left[ \frac{\partial }{\partial \theta _{q}}\left( \dot{%
\tau}\dot{g}\dot{G}_{\cdot h}^{\prime }\right) \right] & 0 \\ 
0 & E^{\ast }\left[ \frac{\partial }{\partial \theta _{q}}\left( \dot{\tau}%
\dot{G}_{\cdot h}\dot{g}^{\prime }\right) \right] & 0 & 0 \\ 
0 & 0 & 0 & 0%
\end{array}%
\right] +$

$\qquad +\left[ 
\begin{array}{cccc}
0 & 0 & E^{\ast }\left[ \frac{\partial }{\partial \theta _{h}}\dot{G}_{\cdot
h}^{\prime }\right] & 0 \\ 
0 & 0 & 0 & 0 \\ 
E^{\ast }\left[ \frac{\partial }{\partial \theta _{h}}\dot{G}_{\cdot h}%
\right] & 0 & -E^{\ast }\left[ \frac{\partial }{\partial \theta _{q}}\left( 
\dot{\tau}\dot{g}\dot{G}_{\cdot h}^{\prime }+\dot{\tau}\dot{G}_{\cdot h}\dot{%
g}^{\prime }\right) \right] & 0 \\ 
0 & 0 & 0 & 0%
\end{array}%
\right] $

\qquad $=\left[ 
\begin{array}{cccc}
0 & 0 & E^{\ast }\left[ \frac{\partial ^{2}\dot{g}^{\prime }}{\partial
\theta _{q}\partial \theta _{h}}\right] & 0 \\ 
0 & -E^{\ast }\left[ \frac{\partial ^{2}\left( \dot{g}\dot{g}^{\prime
}\right) }{\partial \theta _{q}\partial \theta _{h}}\right] & E^{\ast } 
\left[ \frac{\partial ^{2}\left( \dot{g}\dot{g}^{\prime }\right) }{\partial
\theta _{q}\partial \theta _{h}}\right] & 0 \\ 
E^{\ast }\left[ \frac{\partial ^{2}\dot{g}}{\partial \theta _{q}\partial
\theta _{h}}\right] & E^{\ast }\left[ \frac{\partial ^{2}\left( \dot{g}\dot{g%
}^{\prime }\right) }{\partial \theta _{q}\partial \theta _{h}}\right] & 
-E^{\ast }\left[ \frac{\partial ^{2}\left( \dot{g}\dot{g}^{\prime }\right) }{%
\partial \theta _{q}\partial \theta _{h}}\right] & 0 \\ 
0 & 0 & 0 & 0%
\end{array}%
\right] $

It follows that

$\bar{\Psi}_{j}\left( \Phi _{l_{\tau }+h,jk\left( l_{\theta }+q\right)
}^{ETEL}-\Phi _{l_{\tau }+h,jk\left( l_{\theta }+q\right) }^{EL}\right) \bar{%
\Psi}_{k}\xi _{jk}=\mathbf{0}$

$\bar{\Psi}_{j}\left( \Phi _{l_{\kappa }+h,jk\left( l_{\theta }+q\right)
}^{ETEL}-\Phi _{l_{\kappa }+h,jk\left( l_{\theta }+q\right) }^{EL}\right) 
\bar{\Psi}_{k}\xi _{jk}=\mathbf{0}$

$\bar{\Psi}_{j}\left( \Phi _{l_{\lambda }+h,jk\left( l_{\theta }+q\right)
}^{ETEL}-\Phi _{l_{\lambda }+h,jk\left( l_{\theta }+q\right) }^{EL}\right) 
\bar{\Psi}_{k}\xi _{jk}=\bar{\Psi}_{j}\left( \frac{\partial ^{3}}{\partial
\theta _{q}\partial \beta _{j}\partial \beta _{k}}\left( \left( \tau
-1\right) \dot{g}_{h}+\left( \dot{\tau}-\dot{\varepsilon}\right) \dot{g}\dot{%
g}^{\prime }\kappa \right) \right) \xi _{jk}\bar{\Psi}_{k}=$

$\qquad =\left[ 
\begin{array}{c}
0 \\ 
-P\bar{g} \\ 
-P\bar{g} \\ 
-H\bar{g}%
\end{array}%
\right] ^{\prime }\left[ 
\begin{array}{cccc}
0 & 0 & 0 & \left( \frac{3}{2}\right) E^{\ast }\left[ \frac{\partial ^{2}%
\dot{g}_{h}}{\partial \theta ^{\prime }\partial \theta _{q}}\right] \\ 
0 & -2\left( 3\right) E^{\ast }\left[ \frac{\partial \left( \dot{g}_{h}\dot{g%
}\dot{g}^{\prime }\right) }{\partial \theta _{q}}\right] & \left( 3\right)
E^{\ast }\left[ \frac{\partial \left( \dot{g}_{h}\dot{g}\dot{g}^{\prime
}\right) }{\partial \theta _{q}}\right] & 0 \\ 
0 & \left( 3\right) E^{\ast }\left[ \frac{\partial \left( \dot{g}_{h}\dot{g}%
\dot{g}^{\prime }\right) }{\partial \theta _{q}}\right] & 0 & 0 \\ 
\left( \frac{3}{2}\right) E^{\ast }\left[ \frac{\partial ^{2}\dot{g}%
_{h}^{\prime }}{\partial \theta \partial \theta _{q}}\right] & 0 & 0 & 0%
\end{array}%
\right] \left[ 
\begin{array}{c}
0 \\ 
-P\bar{g} \\ 
-P\bar{g} \\ 
-H\bar{g}%
\end{array}%
\right] $

$\qquad =0$, where the parenthesized terms $\left( 3\right) $ and $\left( 
\frac{3}{2}\right) $ arise from $\xi _{jk}$

$\bar{\Psi}_{j}\left( \Phi _{l_{\lambda }+h,jk\left( l_{\theta }+q\right)
}^{ETEL}-\Phi _{l_{\lambda }+h,jk\left( l_{\theta }+q\right) }^{EL}\right) 
\bar{\Psi}_{k}\xi _{jkq}=\bar{\Psi}_{j}\left( \frac{\partial ^{3}}{\partial
\theta _{q}\partial \beta \partial \beta ^{\prime }}\left( \left( \dot{\tau}-%
\dot{\varepsilon}\right) \dot{G}_{\cdot h}^{\prime }\kappa +\dot{\tau}\dot{G}%
_{\cdot h}^{\prime }\lambda \dot{g}^{\prime }\kappa +\left( \tau -\dot{\tau}%
\right) \dot{G}_{\cdot h}^{\prime }\lambda \right) \right) \xi _{jk}\bar{\Psi%
}_{k}=$

$\qquad =\left[ 
\begin{array}{cccc}
0 & -\bar{g}^{\prime }P & -\bar{g}^{\prime }P & -\bar{g}^{\prime }H^{\prime }%
\end{array}%
\right] \left[ 
\begin{array}{cccc}
0 & 0 & \left( 3\right) E^{\ast }\left[ \frac{\partial ^{2}\dot{g}^{\prime }%
}{\partial \theta _{q}\partial \theta _{h}}\right] & 0 \\ 
0 & -\left( 3\right) E^{\ast }\left[ \frac{\partial ^{2}\left( \dot{g}\dot{g}%
^{\prime }\right) }{\partial \theta _{q}\partial \theta _{h}}\right] & 
\left( 3\right) E^{\ast }\left[ \frac{\partial ^{2}\left( \dot{g}\dot{g}%
^{\prime }\right) }{\partial \theta _{q}\partial \theta _{h}}\right] & 0 \\ 
\left( 3\right) E^{\ast }\left[ \frac{\partial ^{2}\dot{g}}{\partial \theta
_{q}\partial \theta _{h}}\right] & \left( 3\right) E^{\ast }\left[ \frac{%
\partial ^{2}\left( \dot{g}\dot{g}^{\prime }\right) }{\partial \theta
_{q}\partial \theta _{h}}\right] & -\left( 3\right) E^{\ast }\left[ \frac{%
\partial ^{2}\left( \dot{g}\dot{g}^{\prime }\right) }{\partial \theta
_{q}\partial \theta _{h}}\right] & 0 \\ 
0 & 0 & 0 & 0%
\end{array}%
\right] \left[ 
\begin{array}{c}
0 \\ 
-P\bar{g} \\ 
-P\bar{g} \\ 
-H\bar{g}%
\end{array}%
\right] $

$\qquad =0$

This implies that $\Phi _{l_{\theta }+l,h}^{-1}\left( \Phi
_{h,jkq}^{ETEL}-\Phi _{h,jkq}^{EL}\right) \bar{\Psi}_{j}\bar{\Psi}_{k}\xi
_{jk}=0$ and thus that

$\qquad E\left[ \left( \Psi _{l_{\theta }+l,jkh}^{ETEL}-\Psi _{l_{\theta
}+l,jkh}^{EL}\right) \bar{\Psi}_{j}\bar{\Psi}_{k}\bar{\Psi}_{h}\bar{\Psi}%
_{l_{\theta }+m}\right] =0+o\left( n^{-2}\right) $.

\section{Conclusion}

$\bar{\Psi}_{\cdot }^{ETEL}=\bar{\Psi}_{\cdot }^{EL}$

$\bar{Q}_{\cdot }^{ETEL}=\bar{Q}_{\cdot }^{EL}$

$\bar{R}_{\cdot }^{ETEL}\not=\bar{R}_{\cdot }^{EL}$ but $E\left[ \left( \bar{%
R}_{l_{\theta }+l}^{ETEL}-\bar{R}_{l_{\theta }+l}^{EL}\right) \bar{\Psi}%
_{l_{\theta }+m}\right] =o\left( n^{-2}\right) $.

In particular, the four terms entering $\left( \bar{R}_{l_{\theta
}+l}^{ETEL}-\bar{R}_{l_{\theta }+l}^{EL}\right) $ have the following
properties

$\left( \bar{\Psi}_{l_{\theta }+l,j}^{ETEL}-\bar{\Psi}_{l_{\theta
}+l,j}^{EL}\right) \bar{Q}_{j}=\frac{1}{2}H_{l\cdot }\bar{g}\bar{g}^{\prime
}P\bar{g}$ (see Section \ref{Rterm1}).

$\left( \Psi _{l_{\theta }+l,jk}^{ETEL}-\Psi _{l_{\theta }+l,jk}^{EL}\right) 
\bar{\Psi}_{k}\bar{Q}_{j}=-\frac{1}{2}H_{l\cdot }\bar{g}\bar{g}^{\prime }P%
\bar{g}-\Xi _{7,l}$ where $\Xi _{7,l}$ is such that $E\left[ \Xi _{7,l}\bar{%
\Psi}_{l_{\theta }+m}\right] =o\left( n^{-2}\right) $

\qquad (see Section \ref{Rterm2}).

$\left( \bar{\Psi}_{l,jk}^{ETEL}-\bar{\Psi}_{l,jk}^{EL}\right) \bar{\Psi}_{j}%
\bar{\Psi}_{k}=0$ (see Section \ref{Rterm3}).

$\left( \Psi _{l_{\theta }+l,jkh}^{ETEL}-\Psi _{l_{\theta
}+l,jkh}^{EL}\right) \bar{\Psi}_{j}\bar{\Psi}_{k}\bar{\Psi}_{h}=\Xi _{8,l}$
such that $E\left[ \Xi _{8,l}\bar{\Psi}_{l_{\theta }+m}^{EL}\right] =o\left(
n^{-2}\right) $ (see Section \ref{Rterm4}).

\end{document}